\documentclass[a4paper,11pt]{article}

\usepackage{amsmath, amssymb, amsfonts, xypic , amsthm, amscd, enumerate, vmargin}
\usepackage[latin1]{inputenc}
\usepackage[french]{babel}

\setmarginsrb{2.4cm}{3.5cm}{2.4cm}{3.3cm}{0cm}{0cm}{0cm}{2.5cm}
 
\theoremstyle{plain}
\newtheorem{teo}{Th\'eor\`eme}[section]
\newtheorem{lem}[teo]{Lemme}
\newtheorem{prop}[teo]{Proposition}
\newtheorem{cor}[teo]{Corollaire}

\theoremstyle{definition}
\newtheorem{defi}[teo]{Definition}

\newtheorem{rem}[teo]{Remarque}

\newcommand{\Hom}{\operatorname{Hom}}

\newcommand{\Spec}{\operatorname{Spec}}

\newcommand{\Sym}{\operatorname{Sym}}

\newcommand{\Cbb}{{\mathbb C}}

\newcommand{\Zbb}{{\mathbb Z}}
\newcommand{\Nbb}{{\mathbb N}}
\newcommand{\Ocal}{{\mathcal O}}
\newcommand{\Bcal}{{\mathcal B}}
\newcommand{\Hcal}{{\mathcal H}}

\newcommand{\Lcal}{{\mathcal L}}
\newcommand{\Qcal}{{\mathcal Q}}
\newcommand{\Ecal}{{\mathcal E}}
\newcommand{\Ncal}{{\mathcal N}}

\newcommand{\Pbb}{{\mathbb P}}
\newcommand{\Abb}{{\mathbb A}}

\newcommand{\Mcal}{{\mathcal M}}

\newcommand{\Fcal}{{\mathcal F}}

\newcommand{\Quot}{\mathcal{Q}uot}

\author{Sébastien Jansou}
\title{Le schéma Quot invariant}
\begin{document}
\date{}

\maketitle

\begin{abstract}
~\\
Etant donné un schéma affine $X$ muni de l'action d'un groupe réductif $G$ et
d'un faisceau $\Mcal$ cohérent et $G$-linéarisé, on construit le  ``schéma Quot invariant'' qui paramètre les quotients de $\Mcal$ dont l'espace des sections globales est somme directe de $G$-modules simples avec des multiplicités finies fixées.

On détermine ensuite le schéma Quot invariant dans une situation simple : $X$ est le cône des vecteurs primitifs d'un $G$-module simple et $\Mcal$ est le faisceau libre sur $X$ engendré par un second $G$-module simple. Le schéma Quot invariant obtenu n'a qu'un seul point, le plus souvent réduit.  Les seuls cas où il n'est pas réduit sont obtenus quand $X$ est le cône des vecteurs primitifs du $\operatorname{Spin}(V)$-module $V$, où $V$ est un
espace vectoriel quadratique de dimension (finie) impaire.

\begin{center} \textbf{Abstract}
\end{center}

Given an affine scheme $X$ with an action of a reductive group $G$ and a $G$-linearized coherent sheaf $\Mcal$, we construct the ``invariant Quot scheme'' that parametrizes the quotients of $\Mcal$ whose space of global sections is a direct sum of simple $G$-modules with fixed finite multiplicities.

Then we determine the invariant Quot scheme in a simple situation, where $X$ is the cone of primitive vectors of a simple $G$-module and $\Mcal$ is the free sheaf on $X$ generated by another simple $G$-module. This invariant Quot scheme has only one point, that is reduced in most of the cases. The only cases where it is not reduced occur when $X$ is the cone of primitive vectors of a quadratic vector space $V$ of odd dimension, under the action of $\operatorname{Spin}(V)$.
\end{abstract}

\section*{Introduction}

Le schéma de Hilbert et le schéma Quot sont des objet fondamentaux en géométrie algébrique. Ils paramètrent respectivement les sous-schémas fermés d'un espace projectif qui admettent un polynôme de Hilbert fixé, et les quotients d'un faisceau cohérent fixé sur un espace projectif qui admettent un polynôme de Hilbert fixé.\\

Haiman et Sturmfels ont obtenu dans \cite{HaS} par des méthodes d'algèbre commutative des objets plus généraux : le \textit{schéma de Hilbert multigradué}, qui paramètre les idéaux homogènes d'une algèbre de polynômes $S$ multigraduée par un groupe abélien qui admettent une ``fonction de Hilbert'' fixée et le \textit{schéma Quot multigradué}, qui paramètre les sous-modules homogènes d'un $S$-module gradué fini qui admettent une fonction de Hilbert fixée.\\

Alexeev et Brion ont construit, à partir du schéma de Hilbert multigradué, le \textit{schéma de Hilbert invariant} : étant donnés un groupe réductif connexe complexe $G$ et une variété $X$ \textit{affine} munie d'une action de $G$, le schéma de Hilbert invariant paramètre les sous-schémas fermés $G$-stables de $X$ dont l'algèbre affine en tant que $G$-module est somme directe de modules simples avec des multiplicités finies fixées. La donnée de ces multiplicités est l'analogue dans cette situation à celle du polynôme de Hilbert.\\

Dans ce travail on vérifie d'abord que, comme dans le cas de la géométrie projective, la construction du schéma de Hilbert invariant se généralise à celle d'un \textit{schéma Quot invariant}, qui paramètre les quotients d'un faisceau $\Mcal$ cohérent $G$-linéarisé fixé sur $X$ par un sous-faisceau $G$-stable tels que l'espace des sections globales du quotient, en tant que $G$-module, soit somme directe de modules simples avec des multiplicités finies fixées. On utilise pour cela le schéma Quot multigradué, qui correspond au cas où le groupe $G$ est un tore.\\

On détermine ensuite une famille ``simple'' de schémas Quot invariants.
Notons $V(\lambda)$ un $G$-module simple de plus grand poids $\lambda$, et
$C_\lambda$ le cône des vecteurs primitifs de $V(\lambda)$, réunion de
l'orbite des vecteurs de plus grand poids et de l'origine. On a montré dans
\cite{hl} que le schéma de Hilbert invariant qui paramètre les déformations
$G$-invariantes de $C_\lambda$ dans $V(\lambda)$ était presque toujours un
point réduit, sauf dans 9 familles de cas où c'est une droite affine. Cela nous
donne une première famille d'exemples de schémas Quot invariants. Ici, on prend comme $G$-variété $X$ le cône $C_\lambda$. C'est le plus petit cône de $V(\lambda)$ stable par $G$. En d'autres termes, son algèbre est la plus petite algèbre graduée engendrée par le $G$-module simple dual de $V(\lambda)$ (on note $\lambda^*$ son plus grand poids), et on a $$\Cbb[C_\lambda] = \bigoplus _ {m \in \Nbb} V(m\lambda^*).$$
On prend comme module $\Mcal$ le module libre $$\Ocal_{C_\lambda} \otimes_\Cbb V(\mu^*)$$  engendré par un module simple $V(\mu^*)$ : l'espace de ses sections globales est $$M=\Cbb[C_\lambda] \otimes_\Cbb V(\mu^*).$$
Les quotients $G$-linéarisés de $\Mcal$ sont les faisceaux dont le module des sections globales est engendré par le module simple $V(\mu^*)$.
On étudie le schéma Quot invariant $\operatorname{Quot}^G(\lambda,\mu)$ des quotients de $M$ qui admettent la décomposition en modules simples $$M/N = \bigoplus _ {m \in \Nbb} V(m\lambda^* + \mu^*).$$
Les multiplicités choisies sont minimales : en effet, si $x$ est un élément primitif de $M/N$ de poids $\mu^*$ et $a$ un élément primitif de $\Cbb[C_\lambda]$ de poids $\lambda^*$, alors l'élément $a^m x $ de $M/N$ est primitif de poids $m\lambda^* + \mu^*$.\\

Un point particulier de
$\operatorname{Quot}^G(\lambda,\mu)$ correspond à la structure de
$\Cbb[C_\lambda]$-module 
$$\bigoplus_{m \in \Nbb} V(m\lambda^* )   \otimes_\Cbb       \bigoplus_{n \in \Nbb} V(n\lambda^* + \mu^*)\longrightarrow \bigoplus_{n \in \Nbb} V(n\lambda^* + \mu^*)$$
donnée sur les composantes homogènes par le produit de Cartan
$$ V(m\lambda^* )   \otimes_\Cbb       V(n\lambda^* + \mu^*) \longrightarrow  V((m+n)\lambda^* + \mu^*).$$
On montrera que le schéma Quot invariant n'admet pas d'autre point
(Proposition \ref{2.6}).
Le plus souvent, le schéma $\operatorname{Quot}^G(\lambda,\mu)$ est un point
réduit. Sinon, c'est le point épais d'ordre 2,  $\Spec( \Cbb[t]/ \langle t^2 \rangle).$ Dans ce cas le module $V(\lambda)$ est le
$\operatorname{Spin}(V)$-module $V$, où $V$ est un espace vectoriel
quadratique de dimension (finie) impaire (théorème \ref{quot}). \\

\noindent \textit{Remerciements :} Je souhaite remercier chaleureusement mon
directeur de thèse Michel Brion pour
son aide.

\section{Construction du schéma Quot invariant}

\subsection{Notations et définition du schéma Quot invariant}

On considère des schémas et des groupes algébriques sur $\Cbb$.
Les références utilisées sont \cite{H} pour la théorie des schémas et \cite{PV} pour celle des groupes algébriques de transformations.

Soit $G$ un groupe réductif connexe. On en choisit un sous-groupe de Borel $B$, et un tore maximal $T$ inclus dans $B$. On considère le radical unipotent $U$ de $B$ : on a : $B=TU$.
Les algèbres de Lie respectives de $G$, $T$ et $U$ sont notées : $\frak{g}$, $\frak{t}$, et $\frak{u} $.
Le système de racines de $G$ relativement à $T$ est noté $R$.
Le choix de $B$ nous en fournit une base $S$, et on a : $R=R_+ \amalg R_-$ où $R_+$ est l'ensemble des racines positives, et $R_-$ celui des racines négatives.

On note $\Lambda$ le groupe des caractères de $T$.
Si $V$ est un $T$-module rationnel (éventuellement de dimension infinie), on note $V=\bigoplus_{ \lambda \in  \Lambda} V_\lambda $ sa décomposition en sous-espaces propres.
Par exemple, l'algèbre de Lie de $G$ admet la décomposition : $$\frak{g}=\displaystyle{  \frak{t}  \oplus \bigoplus_{ \alpha \in R}\frak{g}_\alpha },$$ où chaque $\frak{g}_\alpha$ est de dimension 1. On choisit pour toute racine $ \alpha$ de $G$ un générateur $e_\alpha$ de $\frak{g}_\alpha$.

On a un ordre partiel sur $\Lambda$ : on a $\mu \leq \lambda$ si et seulement si $\lambda - \mu$ est une somme de racines positives.
On note $\Lambda^+$ l'ensemble des éléments de $\Lambda$ qui sont des poids dominants (relativement à la base $S$ du système de racines $R$).

Si $\lambda$ est un poids dominant, on note $V(\lambda)$ le $G$-module dual du $G$-module $$\{ f \in \Cbb[G]~ |~ \forall g \in G,~ \forall b \in B ,~ f(gb)=\lambda(b)f(g)\},$$ où $G$ agit par translations à gauche.
Le $G$-module $V(\lambda)$ est simple, et l'application $\lambda \longmapsto V(\lambda)$ donne une bijection entre les poids dominants de $G$ et les classes d'isomorphisme de $G$-modules simples. Pour l'action de $T$, le poids $\lambda$ est le plus grand poids de $V(\lambda)$. On note $v_\lambda$ le vecteur de $V(\lambda)$ de poids $\lambda$ donné par l'évaluation en l'élément neutre de $G$ : $$v_\lambda  : f \longmapsto f(e).$$

L'algèbre affine du quotient catégorique $G//U$ s'identifie à l'algèbre des invariants $\Cbb[G]^U$. C'est une algèbre graduée par le monoïde des poids dominants : $$\Cbb[G//U]= \bigoplus _ {\lambda \in \Lambda^+} V(\lambda)^*.$$

Si $\lambda$ est un poids dominant, on note $\lambda^*$ le plus grand poids du module dual $V(\lambda)^*$. On a $\lambda^*=-w_0(\lambda)$, où $w_0$ est l'élément le plus long du groupe de Weyl de $G$ relativement au tore maximal $T$.

Si $V$ est un $G$-module rationnel, on note  $V_{(\lambda)}$ sa composante isotypique de type $ \lambda$, c'est-à-dire le sous-module de $V$ somme des sous-modules isomorphes à $ V(\lambda ) $. On a alors la décomposition $V=\bigoplus_{ \lambda \in  \Lambda ^ +} V_{(\lambda)} $. Dans toute décomposition de $V$ en somme directe de modules simples, la multiplicité du module simple $V(\lambda)$ est la dimension de $V^U_\lambda$. Lorsque chacune de ces multiplicités est finie, on dit que le $G$-module $V$ est \textit{à multiplicités finies}.

On sait que les sous-groupes paraboliques de $G$ qui contiennent $B$ sont en bijection avec les parties de $S$. On note $P_I$ le sous-groupe parabolique correspondant à $I \subseteq S$ : les éléments de $I$ sont les racines simples $\alpha$ telles que $- \alpha $ est une racine de $P_I$.

On note $$P_I=L_I U_I$$ la décomposition de Lévi de $P_I$ relativement au tore maximal $T$. Le groupe $L_I$ est le sous-groupe réductif de $G$ qui contient $T$ et dont les racines sont les éléments de $R$ qui sont combinaison linéaire des éléments de $I$. Le groupe $U_I$ est le sous-groupe unipotent de $G$ qui est normalisé par $T$ et dont les racines sont les éléments de $R_+$ qui ne sont pas combinaison linéaire des éléments de $I$.
Si $\lambda$ est un poids dominant, comme le groupe $L_I$ normalise $U_I$, il agit sur l'espace des invariants $$V(\lambda)^{U_I},$$ qui est en fait un $L_I$-module simple de plus grand poids $\lambda$. Il est engendré par les $T$-vecteurs propres de $V(\lambda)$ dont le poids s'écrit $\lambda  - \sum _{\alpha \in I} n_\alpha \alpha$, où les $n_\alpha$ sont des entiers (nécessairement positifs ou nuls).

Soit $\lambda$ un poids dominant, et $P$ un sous-groupe parabolique de $G$ contenant $B$ tel que $\lambda$ se prolonge en un caractère de $P$ (si $P=P_I$, cela signifie que les éléments de $I$ sont orthogonaux à $\lambda$). On note $\pi : G \rightarrow G/P$ la surjection canonique, et $\Lcal _ \lambda$ le faisceau inversible $G$-linéarisé sur $G/P$ qui associe à tout ouvert $\Omega \subseteq G/P$ : $$\Lcal _ \lambda ( \Omega ) := \{ f \in \Ocal _G ( \pi ^ {- 1} ( \Omega) ) \mbox{ $\vert$ } \forall g \in G \mbox{, } \forall p \in P \mbox{, } f(gp)= \lambda ( p ) f (g) \}.$$
L'espace des sections globales de $\Lcal _ \lambda$ est le dual du $G$-module $V(\lambda)$.

Un \textit{$G$-schéma affine} est un schéma affine $X= \Spec A$ de type fini, muni d'une action régulière de $G$. Algébriquement, cela signifie que $A$ est une $\Cbb$-algèbre de type fini sur laquelle $G$ agit par des automorphismes d'algèbre et que pour cette action $A$ est un $G$-module rationnel.
On dit alors que $A$ est une \textit{$G$-algèbre}.  Si $V$ est un $G$-module rationnel de dimension finie, on identifie $V$ au $G$-schéma affine $\Spec (\Sym_\Cbb V^*)$, où $\Sym_\Cbb V^*= \bigoplus _ {n=0}^\infty \operatorname{S}^n V^*$ est l'algèbre symétrique de $ V^*$.

On se fixe un $G$-schéma affine $X= \Spec A$, et un $\Ocal_X$-module $\Mcal$ cohérent et $G$-linéarisé.
La donnée de $\Mcal$ revient à celle de l'espace $M$ de ses sections globales, muni de structures de $A$-module de type fini, et de $G$-module rationnel telles que $$ \forall g \in G \mbox{, }  \forall a  \in A \mbox{, } \forall m \in M  \mbox{, } g ( a m ) = ( g a ) ( g m ).$$ 
On dit alors que $M$ est un \textit{$A$-$G$-module}.
On a un isomorphisme de $A^G$-$G$-modules $$M \cong \bigoplus_{\lambda \in \Lambda^+} \operatorname{Hom}^G(V(\lambda), M) \otimes_\Cbb V(\lambda),$$ en associant à $\sum_\lambda u_\lambda \otimes x_\lambda$ l'élément $\sum_\lambda u_\lambda (x_\lambda)$.

Soit $h : \Lambda ^ + \longrightarrow \Nbb$ une fonction.
On va définir un foncteur contravariant $$\Quot _h ^G (X,  \Mcal) : ( \mbox{Schémas}) ^ \circ \longrightarrow (\mbox{Ensembles}).$$
Soit $S$ un schéma~; on le munit de l'action triviale du groupe $G$.
Notons $\pi$ et $f$ les projections canoniques :

$$\xymatrix{ S \times X   \ar[r]^{~~f} \ar[d]_{\pi} & X \ar[d] \\ S  \ar[r] & \Spec \Cbb } $$

L'image réciproque $f^*  \Mcal$ est un faisceau $G$-linéarisé sur $S \times X $.
Le foncteur $\Quot _h ^G (X,  \Mcal)$ associe au schéma $S$ l'ensemble des sous-$\Ocal_{S \times X}$-modules $\Ncal$ de $f^*\Mcal$ qui sont $G$-stables et tels que l'on ait un isomorphisme de $\Ocal_S$-$G$-modules :
$$ \pi_* ( (f^* \Mcal) / \Ncal) \simeq \bigoplus _ {\lambda \in \Lambda^+} \Fcal_\lambda \otimes_\Cbb V(\lambda),$$ où chaque $\Fcal_\lambda$ est un $\Ocal_S$-module localement libre de rang $h(\lambda)$.
Le quotient $(f^* \Mcal) / \Ncal$ est alors plat sur $S$.

L'objectif des paragraphes \ref{sqm} et \ref{sqi} est d'établir le théorème suivant :

\begin{teo}\label{1.1} Le foncteur $\Quot _h ^G (X,  \Mcal)$ est représenté par un schéma quasi-projectif $\operatorname{Quot} _h ^G (X,  \Mcal)$.
\end{teo}

Le schéma $\operatorname{Quot} _h ^G (X,  \Mcal)$ ainsi défini est appelé le \textit{schéma Quot invariant} des quotients de $\Mcal$ de fonction de Hilbert $h$.

Lorsque $\Mcal$ est le faisceau structural de $X$, on retrouve le schéma de Hilbert invariant.

\subsection{Le schéma Quot multigradué}
\label{sqm}

Comme annoncé dans \cite{HaS} \S6.2, les arguments de la construction du schéma de Hilbert multigradué faite dans cet article se généralisent facilement à la construction d'un \textit{schéma Quot invariant}.
On donne dans cette partie les principales étapes de cette construction.
On commence par rappeler les notations de \cite{HaS}. Comme cette partie traite d'algèbre commutative, on a préféré considérer ici (comme dans \cite{HaS}) des foncteurs covariants de la catégorie des $\Cbb$-algèbres commutatives vers la catégorie des ensembles. (Ce point de vue est équivalent à celui des foncteurs contravariants de la catégorie des schémas vers celle des ensembles.)

Notons $S := \Cbb [x_1, ... , x_n]$ l'algèbre des polynômes à $n$ indéterminées sur $\Cbb$, et $M := \bigoplus _{i=1} ^r S e_i$ un $S$-module libre muni d'une base $\Bcal = (e_1, ... ,e_r)$.
Un \textit{monôme de $M$} est un élément de $M$ de la forme $x^ \alpha e_i$ où $x^ \alpha$ est un monôme de $S$ et $e_i$ un élément de $\Bcal$.

Soit $A$ un groupe abélien. Soit $\deg  : \Nbb^n \longrightarrow A$ un morphisme de monoïdes, et $b_1, ... ,b_n$ des éléments de $A$.
Le \textit{degré} d'un monôme $x^\alpha$ de $S$ (resp.  $x^\alpha e_i$ de $M$) est par définition $\deg \alpha$ (resp. $b_i + \deg \alpha$).
Si $a$ est un élément du groupe $A$, on note $S_a$ (resp. $M_a$) le sous-$\Cbb$-espace vectoriel de $S$ (resp. de $M$) engendré par les monômes de degré $a$.
On obtient ainsi des \textit{multigraduations} de la $\Cbb$-algèbre $S$ et du $S$-module $M$ par le groupe $A$ :
$$ S= \bigoplus_{ a \in A} S_a ~ \mbox{ et } ~ M = \bigoplus_{ a \in A} M_a.$$
Elles vérifient $S_a \cdot S_b \subseteq S_{a+b}$ et $S_a \cdot M_b \subseteq M_{a+b}$.

Le but de cette partie est de paramétrer, une fonction $h : A \longrightarrow \Nbb$ étant donnée, les sous-$S$-modules homogènes $$N = \bigoplus_{ a \in A} N_a \subseteq \bigoplus_{ a \in A} M_a$$ tels que la dimension de $M_a/N_a$ est $h(a)$, pour tout $a \in A$.

Avant de formuler plus précisément le problème, on se place dans le cadre plus général des \textit{$\Cbb$-espaces vectoriels avec opérateurs}.

Un \textit{$\Cbb$-espace vectoriel avec opérateurs} est un $\Cbb$-espace vectoriel $$T = \bigoplus_{ a \in E} T_a$$  gradué par un ensemble $E$ quelconque, et muni d'un ensemble de morphismes de $\Cbb$-espaces vectoriels  $$F= \bigcup _ {a,b \in E} F_{ab}$$ avec $F_{ab} \subseteq \Hom_\Cbb (T_a, T_b)$. On suppose de plus que si $a,b,c$ sont des éléments de $E$, on a $F_{bc} \circ F_{ab} \subseteq F_{ac}$ et que l'application identité sur $T_a$ appartient à $F_{aa}$.\\
Si $D$ est une partie de $E$, on note $T_D$ le $\Cbb$-espace vectoriel gradué $$T_D  := \bigoplus_{ a \in D} T_a, $$ que l'on munit de l'ensemble d'opérateurs $$F_D= \bigcup _ {a,b \in D} F_{ab}.$$

Si $T$ est un $\Cbb$-espace vectoriel gradué muni d'un ensemble $F$ d'opérateurs, et $R$ une $\Cbb$-algèbre commutative, on obtient un $R$-module gradué avec opérateurs par extension des scalaires : le $R$-module gradué est $$R \otimes _ \Cbb T = \bigoplus_{ a \in E} R \otimes _ \Cbb T_a,$$ et on le munit des ensembles d'opérateurs $R$-linéaires $$\widehat{F}_{ab} : R \otimes _\Cbb T_a \longrightarrow R \otimes _\Cbb T_b$$ déduits canoniquement des ensembles d'opérateurs linéaires ${F}_{ab} :  T_a \longrightarrow  T_b$.

Un \textit{sous-$F$-module} de $R \otimes _\Cbb T$ est un sous-$R$-module homogène $$L= \bigoplus_{ a \in E} L_a \subseteq \bigoplus_{ a \in E} R \otimes _\Cbb T_a$$ tel que si $a,b$ sont des éléments de $E$, on a $$\widehat{F}_{ab}(L_a) \subseteq L_b.$$

Soit $h : E \longrightarrow \Nbb$ une fonction.
Pour toute $\Cbb$-algèbre $R$, on note $\Hcal_T^h(R)$ l'ensemble des sous-$F$-modules $L \subseteq R \otimes _ \Cbb T$ tels que le $R$-module $$(R\otimes _ \Cbb T_a) / L_a$$ est localement libre de rang $h(a)$, pour tout $a \in A$.
Si de plus $R \xrightarrow{ ~ \phi ~ } R ^ \prime$ est un morphisme de $\Cbb$-algèbres, $R ^ \prime \otimes_R L$ est un sous-module de $R^\prime \otimes _ R T$ (car $R\otimes _ \Cbb T / L$ est un $R$-module plat), qui est en fait un élément de $\Hcal_T^h(R ^ \prime)$.

On obtient ainsi un foncteur covariant $$\Hcal_T^h : ( \Cbb\mbox{-algèbres}) \longrightarrow (\mbox{Ensembles}).$$

~

On peut maintenant formuler le problème.
La multiplication par les monômes de $S$ munit le $\Cbb$-espace vectoriel gradué $M = \bigoplus_{ a \in A} M_a$ d'opérateurs : les éléments de $F_{ab}$ sont les applications 
$$\begin{array}{rlc}  M_a& \longrightarrow& M_b\\ m &\longmapsto & x^\alpha m
\end{array}$$ pour tout monôme $x^\alpha$ de degré $b-a$.
On remarque qu'ainsi, les sous-$F$-modules de $M$ ne sont autres que les sous-$S$-modules homogènes de $M$. Le but de ce paragraphe est d'établir le théorème suivant, qui définit le \textit{schéma Quot multigradué} $H_M^h$ :

\begin{teo}\label{1.2} Soit $h :A \longrightarrow \Nbb$ une fonction.\\ Le foncteur $\Hcal_M^h$ est représenté par un schéma quasi-projectif $H_M^h$.
\end{teo}

La démonstration se fait en deux étapes. On montre d'abord que pour toute partie finie $D$ du groupe abélien $A$, le foncteur $\Hcal_{M_D}^h$ est représentable par un schéma quasi-projectif $H_{M_D}^h$ (Proposition \ref{1.6}). 
On montre ensuite qu'il existe une partie finie $D$ de $A$ telle que le foncteur $\Hcal_M^h$ est représenté par un sous-schéma fermé de $H_{M_D}^h$ (Lemme \ref{1.7} et Proposition \ref{1.8}).
On commence par montrer un lemme combinatoire, utilisé lors de chacune des deux étapes.

~

Un \textit{sous-module monomial} de $M$ est un sous-$S$-module de $M$ engendré par des monômes de $M$.
Les sous-modules monomiaux de $M$ sont donc ceux de la forme $\bigoplus _{i=1} ^r I_i e_i$, où $I_1, ... ,I_r$ sont des idéaux monomiaux de $S$.

On dit qu'un ensemble $\Ecal$ de sous-modules de $M$ est une \textit{antichaîne} si pour tout couple $(N_1, N_2)$ d'éléments de $\Ecal$, on a $N_1 \not \subseteq N_2$. 

Maclagan a montré (\cite{mac}) que les antichaînes d'idéaux monomiaux d'une algèbre de polynômes sont finies. Le lemme suivant en est une généralisation immédiate :

\begin{lem}\label{1.3} Les antichaînes de sous-modules monomiaux de $M$ sont finies.
\end{lem}

\noindent {\bf Preuve.}
Associons à tout sous-module monomial $N=\bigoplus _{i=1} ^r I_i e_i$ de $M$ l'idéal monomial $$J_N  := \sum _{i=1} ^r I_i y_i + \sum _{i,j} y_i y_j \Cbb [x_1, ... , x_n, y_1, ... ,y_r] $$ de l'algèbre de polynômes $\Cbb [x_1, ... , x_n, y_1, ... ,y_r]$. Pour tous sous-modules monomiaux $N_1$, $N_2$ de $M$, on a $N_1 \subseteq N_2$ si et seulement si $J_{N_1} \subseteq J_{N_2}$.
On associe ainsi à toute antichaîne de sous-modules monomiaux de $M$ une antichaîne d'idéaux monomiaux de $\Cbb [x_1, ... , x_n, y_1, ... ,y_r]$, et on en déduit le lemme.
\hfill $\Box$

~

Si $N$ est un $\Cbb$-espace vectoriel de dimension finie et $r$ un entier tel que $0 \leq r \leq \dim N$, on note $G_N^r$ la grassmannienne des quotients de $N$ de dimension $r$.
Si de plus l'espace vectoriel $N=\bigoplus_{a \in E} N_a$  est gradué par un ensemble fini $E$, et $h :E \longrightarrow \Nbb$ est une fonction, on note $G_N^h$ la grassmannienne des quotients de $N$ par un sous-espace vectoriel homogène $N^\prime$ tel que $\dim N_a / N^\prime _ a = h(a)$ pour tout $a \in A$.
Ce schéma est donc un produit de grassmanniennes :$$G_N^h = \prod _{a \in E} G_{N_a}^{h(a)}.$$
On définit enfin, si de plus $M$ est un sous-espace vectoriel de $N$, la \textit{grassmannienne relative} $G_{N \setminus M}^h$. Il s'agit de l'ouvert de $G_N^h$ qui paramètre les quotients $N/N^\prime$ de $N$ qui sont \textit{engendrés par $M$,} c'est-à-dire tels que le morphisme canonique $M \longrightarrow N/N^\prime$ soit surjectif (on renvoie à \cite{HaS} Proposition 2.11 pour plus de détails).

On rappelle ici les deux théorèmes suivants, établis dans \cite{HaS} (theorems 2.2, 2.3) :
\begin{teo}\label{1.4}Soit $(T,F)$ un $\Cbb$-espace vectoriel avec opérateurs dont l'ensemble $E$ des degrés est \textit{fini}. Soit $h :E \longrightarrow \Nbb$ une fonction.
Soit $M \subseteq N \subseteq T$ deux sous-$\Cbb$-espaces vectoriels homogènes de $T$.
Soit $G \subseteq F$ un sous-ensemble.
Supposons
\begin{enumerate}[(1)]
 \item $N$ est un $\Cbb$-espace vectoriel de dimension finie.
 \item $N$ engendre le $F$-module $T$.
 \item Pour tout surcorps $K$ de $\Cbb$, et tout élément $L$ de $\Hcal_T^h(K)$, l'application naturelle\\
 $K \otimes_\Cbb M \longrightarrow K \otimes_\Cbb T/L$ est surjective. 
 \item $G$ engendre $F$ comme catégorie, et $G.M \subseteq N$.
\end{enumerate}
Alors le foncteur $\Hcal_T^h$ est représenté par un sous-schéma fermé de la grassmannienne relative $G_{N \setminus M}^h$, donc par un schéma quasi-projectif.
\end{teo}

\begin{teo}\label{1.5}Soit $(T,F)$ un $\Cbb$-espace vectoriel avec opérateurs, et $h :E \longrightarrow \Nbb$ une fonction. Soit $D$ une partie de $E$ telle que $\Hcal_{T_D}^h$ est représenté par un schéma $H_{T_D}^h$. Supposons que pour tout $a \in E$ :
\begin{enumerate}[(1)]
 \item Il existe une partie $G$ finie de $\bigcup _ { b \in D} F_{ba}$ telle que le $\Cbb$-espace vectoriel $T_a / \sum _{ b \in D} G_{ba}(T_b)$ est de dimension finie.
 \item Pour tout surcorps $K$ de $\Cbb$, et tout élément $L_D$ de $\Hcal_{T_D}^h(K)$, le sous-$F$-module $L^\prime$ de $K \otimes _ \Cbb T$ engendré par $L_D$ vérifie $$ \dim_K ( K\otimes _ \Cbb T_a / L^\prime_a) \leq h(a).$$
\end{enumerate}
Alors $\Hcal_{T}^h$ est représenté par un sous-schéma fermé de $H_{T_D}^h$.
\end{teo}

On obtient la proposition suivante en appliquant le théorème \ref{1.4} à l'aide du lemme \ref{1.3}. La démonstration est analogue à la première partie de \cite{HaS} Proof of theorem 1.1, p 742.
\begin{prop}\label{1.6} Soit $D$ une partie \textit{finie} du groupe abélien $A$.
Le foncteur $\Hcal_{M_D}^h$ est représenté par un schéma quasi-projectif $H_{M_D}^h$.

\end{prop}

Si $N$ est un sous-espace vectoriel homogène de $M$, on note  $h_N(a)$  la dimension (éventuelle\-ment infinie) du quotient $M_a/N_a$,  pour tout $a \in A$.

On déduit le lemme suivant du lemme \ref{1.3}. La démonstration est analogue à celle de \cite{HaS} Proposition 3.2.
\begin{lem}\label{1.7} Il existe une partie finie $D$ de $A$ telle que

\begin{enumerate}[(1)]
 \item Tous les sous-modules monomiaux $N$ tels que $h_N=h$ sont engendrés par leurs éléments homogènes de degré appartenant à $D$.
 \item Si $N$ est un sous-module monomial de $M$ engendré par ses éléments homogènes de degré appartenant à $D$ et tel que $h_N|_D = h|_D$, alors $h_N = h$.
\end{enumerate}
\end{lem}

On obtient enfin la proposition suivante en appliquant le théorème \ref{1.5}. La démonstration est analogue à la seconde partie de \cite{HaS} Proof of theorem 1.1, p 742.

 \begin{prop}\label{1.8}Soit $D$ une partie de $A$ donnée par le lemme \ref{1.7} (en particulier, $D$ est finie).
Alors le foncteur $\Hcal_{M}^h$ est représenté par un sous-schéma fermé de $H_{M_D}^h$, donc par un schéma quasi-projectif.
\end{prop}

Le théorème \ref{1.2} est donc démontré.
\subsection{Fin de la construction}
\label{sqi}

Dans cette partie, on donne la construction du schéma Quot invariant (à partir du schéma Quot multigradué), parfaitement analogue à celle du schéma de Hilbert invariant d' Alexeev-Brion. Comme elle ne présente aucune difficulté nouvelle, on s'est contenté de donner les principales étapes, sans preuves complètes.

Traitons d'abord le cas où le groupe $G$ est un tore : on a $G=T$.
Soit $Y$ un $T$-schéma affine, et $\Mcal$ un faisceau cohérent $T$-linéarisé sur $Y$. On note $M$ l'espace des sections globales de $\Mcal$.

Soit $E$ un $T$-module de dimension finie tel que $Y$ s'identifie (en tant que $T$-schéma) à un sous-schéma fermé $T$-stable de $E$. 
Soit $(e_1,...,e_r)$ un système de générateurs fini du $A$-module $M$ formé de vecteurs propres pour l'action de $T$.
On associe à ce système de générateurs une surjection de $\Ocal_E$-modules $T$-linéarisés $$\widetilde{\Mcal} :=\bigoplus_{i=1}^r \Ocal_E e_i     \twoheadrightarrow            \Mcal.$$

Le théorème \ref{1.2} nous donne que le foncteur $\Quot _h ^T (E,  \widetilde{\Mcal})$ est représenté par un schéma quasi-projectif $\operatorname{Quot} _h ^T (E,  \widetilde{\Mcal})$.
 Le lemme suivant correspond au théorème \ref{1.1} dans le cas où le groupe réductif $G$ est un tore.
Sa démonstration est analogue à celle de \cite{AB} lemma 1.6.

\begin{lem}\label{1.9} Le foncteur $\Quot _h ^T (Y,  \Mcal)$ est représenté par un sous-schéma fermé de\\ $\operatorname{Quot} _h ^T (E,  \widetilde{\Mcal})$, donc par un schéma quasi-projectif.
\end{lem}

Traitons maintenant le cas général. On garde les notations du théorème \ref{1.1}.

On note $X//U$ le quotient catégorique du $G$-schéma affine $X= \Spec A$ par le sous-groupe unipotent maximal $U$ de $G$ : $$X//U  := \Spec A^U.$$
(On rappelle que $A^U$ est une $\Cbb$-algèbre de type fini, selon \cite{Gross}, Thm 9.4.)\\
Le schéma affine $X//U$ est muni d'une action du tore $T$.

Notons  $\Mcal ^ U$ le faisceau $T$-linéarisé sur $X//U$ des $U$-invariants du faisceau $\Mcal$. C'est un faisceau cohérent (en effet, montrons que l'espace de ses sections globales $M^U$ est un $A^U$-module de type fini. Comme $M$ est un $A$-module de type fini, son algèbre symétrique  $\Sym_A(M)$ est une $A$-algèbre graduée de type fini, donc une $\Cbb$-algèbre graduée de type fini. L'algèbre de ses $U$-invariants $$\Sym_A(M)^U = A^U \oplus M^U \oplus (\operatorname{S}^2M)^U \oplus ...$$ est donc aussi une $\Cbb$-algèbre graduée de type fini : en particulier, sa composante homogène de degré $1$ est un $A^U$-module de type fini.)

Le foncteur $\Quot _h ^G (X,  \Mcal)$ peut être vu comme un sous-foncteur de $\Quot _h ^T (X//U,  \Mcal^U)$ (on prolonge la fonction $h$ à $\Lambda$ en posant $h=0$ sur $\Lambda \setminus \Lambda^+$).

On a en effet un morphisme fonctoriel $\phi$ donné pour tout schéma $S$ par 
$$\begin{array}{ccc}   \Quot _h ^G (X,  \Mcal)(S) &\xrightarrow{ \phi(S) }& \Quot _h ^T (X//U,  \Mcal^U)(S) \\ \Ncal &\longmapsto & \Ncal^U
\end{array}$$
et les $\phi(S)$ sont des injections car le seul antécédent de $\Ncal^U$ possible est le $G$-module engendré par $\Ncal^U$.

Selon le lemme \ref{1.9}, le foncteur $\Quot _h ^T (X//U,  \Mcal^U)$ est représenté par un schéma quasi-projectif $\operatorname{Quot} _h ^T (X//U,  \Mcal^U)$.

\begin{prop}\label{1.10}
 Le sous-foncteur $\Quot _h ^G (X,  \Mcal) \hookrightarrow \Quot _h ^T (X//U,  \Mcal^U)$ est représenté par un sous-schéma fermé de $\operatorname{Quot}_h ^T (X//U,  \Mcal^U)$, donc par un schéma quasi-projectif.
\end{prop}

La démonstration est analogue à celle du Thm 1.7 de \cite{AB}.

Le théorème \ref{1.1} est donc démontré.

\subsection{Premières propriétés du schéma Quot invariant}

Dans ce paragraphe, on note toujours $X$ un $G$-schéma affine, $\Mcal$ un faisceau cohérent $G$-linéarisé sur $X$ dont on note $M$ l'espace des sections globales, et $h :\Lambda^+ \longrightarrow \Nbb$ une fonction.

La proposition suivante décrit l'espace tangent au schéma Quot invariant en un point fermé. On donne sa démonstration, analogue à celle de \cite{AB} Proposition 1.13, pour expliciter l'isomorphisme canonique.
\begin{prop}\label{1.11}
Soit $z$ un point fermé du schéma $\operatorname{Quot} _h ^G (X,  \Mcal) $, c'est-à-dire un sous-module $\Ncal \subseteq \Mcal$ stable par $G$ et tel que, en notant $N$ l'espace des sections globales de $\Ncal$, on ait un isomorphisme de $G$-modules $$M/N \simeq \bigoplus _ {\lambda \in \Lambda^+} {h(\lambda)}V(\lambda).$$
L'espace tangent de Zariski au schéma Quot invariant est canoniquement isomorphe à l'espace des morphismes de $A$-$G$-modules de $N$ dans $M/N$ : $$T_z \operatorname{Quot} _h ^G (X,  \Mcal) \cong \Hom_A^G (N, M/N).$$
\end{prop}

\noindent {\bf Preuve.}

Notons $\epsilon$ la classe de $t$ dans l'algèbre $\Cbb[t]/\langle t^2 \rangle$.

L'espace tangent en $z$ est l'ensemble des morphismes de $\Spec \Cbb[\epsilon]$ dans $\operatorname{Quot} _h ^G (X,  \Mcal)$ dont la restriction à $\Spec \Cbb$ (vu comme un sous-schéma fermé de $\Spec \Cbb[\epsilon]$) correspond au point $z$.

En d'autres termes, c'est l'ensemble des sous-$A[\epsilon]$-$G$-modules $$L \subseteq \Cbb[\epsilon] \otimes_\Cbb M = M \oplus \epsilon M$$ tels qu'on ait l'identification $$\Cbb \otimes_{\Cbb[\epsilon]} L \cong N$$ et que le quotient $$ (\Cbb[\epsilon] \otimes_\Cbb M) / L $$ soit un $\Cbb[\epsilon]$-module plat.

Soit un tel sous-module $L$.

Précisons d'abord la première condition. On rappelle qu'on a un plongement (grâce à la seconde condition) $$\Cbb \otimes_{\Cbb[\epsilon]} L \hookrightarrow \Cbb \otimes_{\Cbb[\epsilon]} (\Cbb[\epsilon] \otimes_\Cbb M) \cong (M \oplus \epsilon M)/\epsilon M \cong M.$$
La première condition dit que l'image de ce plongement est $N$.
Autrement dit, la projection de $L \subseteq M \oplus \epsilon M$ sur $M$ est $N$, c'est-à-dire
\begin{equation}\label{pullback}
(L+\epsilon M) \cap M = N.      
\end{equation}
Notons que l'on a donc (en multipliant (\ref{pullback}) par $\epsilon$) $\epsilon L = \epsilon N$.

On utilise maintenant la seconde condition. Notons $v$ un élément de $\Cbb[\epsilon] \otimes_\Cbb M$, et $\overline{v}$ sa classe dans le quotient $ (\Cbb[\epsilon] \otimes_\Cbb M) / L $.

La seconde condition signifie que si $\epsilon \overline{v} =0$, alors $  \overline{v}$ appartient à $\epsilon((\Cbb[\epsilon] \otimes_\Cbb M) / L)$.

Autrement dit, si $\epsilon v $ appartient à $L$, alors $v$ appartient à $\epsilon M + L$.

Donc si $\epsilon v $ appartient à $L$, alors $\epsilon v$ appartient à $\epsilon (\epsilon M + L) =\epsilon L= \epsilon N$.

D'où $ L \cap \epsilon M \subseteq \epsilon N$, et comme l'inclusion réciproque est toujours vraie, on a
\begin{equation}\label{platitude}
 L \cap \epsilon M = \epsilon N.
\end{equation}

On peut maintenant conclure. 
Pour tout élément $n$ de $N$, il existe un unique élément $\phi(n)$ de $M/N$ (on voit cet élément comme une partie de $M$) tel que $$n+\epsilon \phi(n) \subseteq L$$ (l'unicité découle de (\ref{platitude}) et l'existence de (\ref{pullback})).
On a alors
 \begin{equation}\label{formule}L= \bigcup_{n \in N} (n + \epsilon \phi(n)).\end{equation}
Comme $L$ est un $A$-$G$-module, l'application $\phi$ est un morphisme de $A$-$G$-modules.

Réciproquement, tout morphisme de $A$-$G$-modules $\phi: N \longrightarrow M/N$ définit bien via l'expression (\ref{formule}) un morphisme de $\Spec \Cbb[\epsilon]$ dans $\operatorname{Quot} _h ^G (X,  \Mcal)$. 
\hfill $\Box$\\

La proposition suivante est une généralisation de \cite{HaS} Corollary 1.2.

\begin{prop}\label{1.12} Supposons que le $G$-module $M$ est à multiplicités finies. Alors le schéma $\operatorname{Quot} _h ^G (X,  \Mcal)$ est projectif.
\end{prop}

\noindent {\bf Preuve.}
Par construction, le schéma $\operatorname{Quot} _h ^G (X,  \Mcal)$ est quasi-projectif.
Pour montrer qu'il est projectif, il suffit donc de montrer qu'il est propre sur $\Cbb$.
Pour cela, on utilise le critère valuatif de propreté (voir \cite{H} Theorem II.4.7).
Soit $R$ un anneau de valuation discrète, et $K$ son corps de fractions.
Il s'agit de montrer que tout morphisme $\Spec K \xrightarrow{\phi} \operatorname{Quot} _h ^G (X,  \Mcal)$ se prolonge en un morphisme $\Spec R \xrightarrow{\widetilde{\phi}} \operatorname{Quot} _h ^G (X,  \Mcal)$. 

Un tel morphisme $\phi$ revient à un sous-$K \otimes_\Cbb A$-module $G$-stable $$L \subseteq K \otimes_\Cbb M$$ tel que pour tout poids dominant $\lambda$, le $K$-espace vectoriel $(K \otimes_\Cbb M / L ) ^U _ \lambda$ est de dimension $h(\lambda)$.

On considère $R \otimes_\Cbb M$ comme un sous-$R\otimes_\Cbb A$-module $G$-stable de $K \otimes_\Cbb M$.
L'espace $$P :=L \cap (R\otimes_\Cbb M)$$ est un sous-$R\otimes_\Cbb A$-module $G$-stable de $R\otimes_\Cbb M$.

De plus, comme le $K$-espace vectoriel $K\otimes_\Cbb M_\lambda^U$ est de dimension finie, le $R$-module $$(R \otimes_\Cbb M / P ) ^U _ \lambda = R \otimes_\Cbb M^U_\lambda /( L_\lambda^U\cap(R\otimes_\Cbb M_\lambda^U) )$$ est lui aussi libre de rang $h(\lambda)$.

Le sous-module $P\subseteq R \otimes_\Cbb M$ correspond donc à un morphisme $\Spec R \xrightarrow{\widetilde{\phi}} \operatorname{Quot} _h ^G (X,  \Mcal)$. Enfin, la restriction de $\widetilde{\phi}$ est bien $\phi$, car on a $K \otimes _ R P = L$ (en effet, par définition de $P$, on a $K \otimes _ R P \subseteq L$, et tout élément de $L$ est égal, à un scalaire appartenant à $K$ près, à un élément de $P$).
\hfill $\Box$

~

Lorsque le groupe $G$ est trivial, le seul poids de $G$ est le poids nul, et la donnée d'une fonction de Hilbert revient donc à celle d'un entier $n \in \Nbb$. On note alors respectivement $\operatorname{Hilb} _n  (X)$ et $ \operatorname{Quot} _n (X,  \Mcal)$ le schéma de Hilbert invariant et le schéma Quot invariant.

Le schéma de Hilbert invariant $\operatorname{Hilb} _n  (X)$ n'est autre que le schéma de Hilbert des sous-schémas de longueur $n$ de $X$ (c'est-à-dire le schéma de Hilbert de $n$ points sur $X$, qui est défini dès que $X$ est un schéma quasi-projectif).

On a naturellement un morphisme fonctoriel $$\Quot _h ^G (X,  \Mcal) \longrightarrow \Quot _{h(0)} (  X//G,  \Mcal^G) :$$ avec les notations du \S1.1, il associe à tout élément $\Ncal \subseteq f^* \Mcal$ de $\Quot _h ^G (X,  \Mcal)(S)$ le sous-$\Ocal_{X//G}$-module $\Ncal^G$ de $(f^* \Mcal)^G$.
On a donc un morphisme naturel de schémas $$\gamma  : \operatorname{Quot} _h ^G (X,  \Mcal) \longrightarrow \operatorname{Quot} _{h(0)} (  X//G,  \Mcal^G).$$

Dans le cas du schéma de Hilbert invariant, c'est-à-dire si $\Mcal$ est le faisceau structural de $X$, ce morphisme associe à tout fermé $G$-stable $Y \subseteq X$ de fonction de Hilbert $h$ le fermé $Y//G$ (qui est en fait fini) de $X//G$. Ce morphisme est donc un analogue au morphisme de Hilbert-Chow de Nakamura (\cite{Na} \S 2.1). Signalons cependant qu'il ne généralise pas le ``morphisme de Chow'' défini par Haiman et Sturmfels pour le cas du schéma de Hilbert torique (\cite{HaS} \S5).

\begin{prop}
Le morphisme $\gamma  : \operatorname{Hilb} _h ^G (X) \longrightarrow \operatorname{Hilb} _{h(0)} (  X//G)$ est projectif.
\end{prop} 

\noindent {\bf Preuve.}
Comme pour la proposition précédente, il suffit de montrer que ce morphisme est propre.
On utilise à nouveau le critère valuatif de propreté : soit $R$ un anneau de valuation discrète, et $K$ son corps de fractions. 

Soient deux morphismes $\phi$ et $\psi$ tels qu'on ait un diagramme commutatif :

$$\xymatrix{  \operatorname{Spec}(K)  \ar[r]^{\phi} \ar[d]_{} & \operatorname{Hilb} _h ^G (X) \ar[d] \\ \operatorname{Spec}(R)  \ar[r]^{\psi~~~~~} & \operatorname{Hilb} _{h(0)} (  X//G) } $$

Il faut montrer que le morphisme $\phi$ se prolonge en un morphisme $\operatorname{Spec}(R)  \xrightarrow{\widetilde{\phi}} \operatorname{Hilb}_h ^G (X) $. 

Le morphisme $\phi$ correspond à un idéal $G$-stable $I \subseteq K \otimes_\Cbb A$ tel que pour tout poids dominant $\lambda$, le $K$-espace vectoriel $(K \otimes_\Cbb A / I ) ^U _ \lambda$ est de dimension $h(\lambda)$.

Le morphisme $\psi$ correspond à un idéal $G$-stable $J \subseteq R \otimes_\Cbb A^G$ tel que le $R$-module $(R \otimes_\Cbb A^G) / J $ est de dimension $h(0)$.

Enfin, la commutativité du diagramme signifie que $$K \otimes _ R J = I^G.$$

Comme précédemment, on considère l'idéal $$J ^\prime  := I \cap (R \otimes _\Cbb A)$$ de $R \otimes _\Cbb A$. Il est stable par $G$.

Montrons que pour tout poids dominant $\lambda$, le $R$-module $$(R \otimes_\Cbb A / {J ^\prime} ) ^U _ \lambda$$ est libre de rang $h(\lambda)$.

On remarque que c'est un module sans torsion, donc plat (car $R$ est principal).
Montrons que c'est un module de type fini.

Pour cela, il suffit de montrer que $$(R \otimes_\Cbb A^U _ \lambda) / (J \cdot (R \otimes_\Cbb A^U _ \lambda))$$ est un $R$-module de type fini, car $J \cdot (R \otimes_\Cbb A^U _ \lambda)$ est inclus dans ${J ^\prime}^U _ \lambda$.

On sait (voir par exemple \cite{AB} Lemma 1.2) que $R\otimes _\Cbb A^U _ \lambda$ est un $R \otimes_\Cbb A^G$-module de type fini. 

Donc le quotient $(R\otimes _\Cbb A^U _ \lambda)/ (J \cdot (R \otimes_\Cbb A^U _ \lambda))$ est un $(R \otimes_\Cbb A^G)/J$-module de type fini, donc un $R$-module de type fini (car $(R \otimes_\Cbb A^G)/J$ est un $R$-module de type fini).

Ainsi, le $R$-module $(R \otimes_\Cbb A / {J ^\prime} ) ^U _ \lambda$ est plat de type fini : il est donc libre (car $R$ est local).

Enfin, on a (comme précédemment) $K \otimes _R J^\prime = I$, donc le rang de $(R \otimes_\Cbb A / {J ^\prime} ) ^U _ \lambda$ est $h(\lambda)$.

 L'idéal $J^\prime$ correspond donc à un morphisme ${\widetilde{\phi}} :\operatorname{Spec}(R)  \longrightarrow \operatorname{Hilb}_h ^G (X) $, dont la restriction à $\Spec K$ est $\phi$. \hfill $\Box$\\

Par contre, dans le cas du schéma Quot invariant, le morphisme $\gamma$ n'est pas nécessairement projectif.
Par exemple, supposons que $G$ est le groupe multiplicatif, et que $X$ est la droite affine $\Abb ^1$ munie de l'action triviale de $G$. Notons $\Cbb_1$ la droite vectorielle où $G$ agit avec le poids $1$, et $h$ la fonction valant $1$ sur le poids $1$ et $0$ ailleurs. Supposons enfin que $\Mcal  := \Ocal_X \otimes_\Cbb \Cbb_1$ est le faisceau libre sur $X$ où $G$ agit avec poids $1$ sur les sections.

Le schéma  $\operatorname{Quot} _h ^G (\Abb ^1,  \Mcal)$ coïncide avec $\operatorname{Hilb} _1 (\Abb ^1)=\Abb ^1$. Le schéma $\operatorname{Quot} _0  (\Abb ^1,  \Mcal^G)$ consiste en un point réduit (le faisceau $\Mcal^G$ est nul). Donc le morphisme $\gamma  : \Abb ^1 = \operatorname{Quot} _h ^G (X,  \Mcal) \longrightarrow \operatorname{Quot} _{h(0)} (  X//G,  \Mcal^G)= \Spec \Cbb$ n'est pas projectif. 

\section{Etude d'une classe de schémas Quot invariants}

Soit $\lambda$ un poids dominant.

On a une action régulière de $G$ sur l'espace $\Pbb (V( \lambda)  ) $ des droites de $V( \lambda)$.
Notons $[ v _ \lambda ] \in \Pbb (V( \lambda)  ) $ la droite engendrée par $v _ \lambda$ et $$P _ \lambda  := G_{[ v _ \lambda ]}$$son stabilisateur dans $G$ : c'est le plus grand sous-groupe parabolique de $G$ qui contient $B$ et tel que $\lambda$ se prolonge en un caractère de $P_\lambda$.
On a donc $P_ \lambda = P_I$, où $I$ est l'ensemble des racines simples qui sont orthogonales à $\lambda$.
On note $$P_\lambda = L_\lambda U_\lambda$$ la décomposition de Lévi de $P_ \lambda $ relativement au tore maximal $T$.  
L'orbite de $[ v _ \lambda ]$ est la seule orbite fermée de $\Pbb (V( \lambda)  ) $ (donc l'unique orbite de plus petite dimension).
L'espace homogène projectif $G / P_ \lambda $ se plonge ainsi dans $\Pbb (V( \lambda)  ) $, et le faisceau inversible très ample associé à ce plongement est en fait $\Lcal _ \lambda$.
Le cône affine au dessus de $G / P_ \lambda $ dans $V( \lambda)$ est le cône $$C_\lambda  := \overline{G \cdot v_\lambda}= G \cdot v_\lambda \cup \{ 0 \}$$ des vecteurs primitifs de $V(\lambda)$.

On note $A(\lambda)$ l'algèbre affine de $C_\lambda$.

Comme le morphisme dominant $$\begin{array}{rlc}  G& \longrightarrow& C_\lambda\\ g &\longmapsto & g \cdot v_\lambda
\end{array}$$ se factorise par $G//U$, l'algèbre $A(\lambda)$ s'identifie à la sous-algèbre de $\Cbb[G//U]$ engendrée par $V(\lambda)^*$ :
$$A(\lambda) \cong \bigoplus _{m=0}^\infty V(m\lambda)^*.$$
C'est une algèbre graduée par $\Nbb$.

Soit $\mu$ un poids dominant.
Notons $Q(\lambda, \mu)$ le sous-$A(\lambda)$-module $G$-stable de
$\Cbb[G//U]$~: $$ Q(\lambda, \mu) :=\bigoplus _{m=0}^\infty V(m\lambda+\mu)^*.$$
C'est un $A(\lambda)$-module gradué par $\Nbb$, engendré par sa composante homogène de degré nul $Q(\lambda, \mu)_0=V(\mu)^*$.

On a donc une surjection de $A(\lambda)$-$G$-modules $$M(\lambda, \mu) := A(\lambda) \otimes_\Cbb V(\mu)^* 
\twoheadrightarrow Q(\lambda, \mu) .$$

Notons $N(\lambda, \mu)$ son noyau.
Les $A(\lambda)$-modules $M(\lambda, \mu)$ et $N(\lambda, \mu)$ sont gradués par $\Nbb$~: la composante homogène de degré $m$ de $M(\lambda, \mu)$ est $$M(\lambda, \mu)_m= V(m\lambda)^* \otimes_\Cbb V(\mu)^*~; $$ celle de $N(\lambda, \mu)$ est le noyau du produit de Cartan $$V(m\lambda)^* \otimes_\Cbb V(\mu)^* \longrightarrow  V(m\lambda+\mu)^* .$$
En particulier, les composantes isotypiques non nulles de $N(\lambda, \mu)_m$ sont toutes de type strictement inférieur à $m \lambda^* + \mu^*$.

On a une suite exacte de $A(\lambda)$-$G$-modules gradués $$0 \longrightarrow N(\lambda, \mu) \longrightarrow M(\lambda, \mu) \longrightarrow Q(\lambda, \mu) \longrightarrow 0.$$
Notons $\Ncal(\lambda, \mu)$, $\Mcal(\lambda, \mu)$, $\Qcal(\lambda, \mu)$ les faisceaux cohérents $G$-linéarisés sur $C_\lambda$ correspondant respectivement à $N(\lambda, \mu)$, $M(\lambda, \mu)$ et $Q(\lambda, \mu)$.
On a donc une suite exacte de faiseaux $G$-linéarisés $$0 \longrightarrow \Ncal(\lambda, \mu) \longrightarrow \Mcal(\lambda, \mu) \longrightarrow \Qcal(\lambda, \mu) \longrightarrow 0.$$

Notons $h_{\lambda, \mu}  : \Lambda ^+ \longrightarrow \Nbb $ la fonction valant $1$ sur les poids de la forme $m \lambda^* + \mu^*$, et $0$ ailleurs. On a vu que le quotient $\Mcal(\lambda, \mu)/\Ncal(\lambda, \mu) \cong \Qcal(\lambda, \mu)$ admet la fonction de Hilbert $h_{\lambda, \mu} $~: il correspond donc à un point fermé du schéma Quot invariant $$\operatorname{Quot}_{h_{\lambda, \mu}}^G(C_\lambda, \Mcal(\lambda, \mu)).$$
 On note désormais ce schéma  $\operatorname{Quot}^G(\lambda, \mu).$

\begin{rem}\label{2.1} Selon la proposition \ref{1.12}, le schéma $\operatorname{Quot}^G(\lambda, \mu)$ est projectif. En effet, comme l'algèbre $A(\lambda)$ est à multiplicités finies, le $G$-module $M(\lambda, \mu) :=A(\lambda) \otimes_\Cbb V(\mu)^*$ est à multiplicités finies.
\end{rem}

Le but de cette partie est de démontrer le théorème suivant.
Les notations utilisées concernant les systèmes de racines sont celles de Bourbaki (\cite{Bou}).

\begin{teo} \label{quot} Le schéma Quot invariant $\operatorname{Quot}^G(\lambda, \mu)$ est un point réduit, sauf si on a (à revêtement fini de $G$ près) $G= \operatorname{Spin}(2n+1) \times H $ pour un groupe réductif
 connexe $H$ et $V(\lambda)=\Cbb^{2n+1}$ et $\mu=(\mu_1, \mu_2)$ avec
 $\langle \mu_1 , \alpha_n^\vee \rangle \geq 1$. 
On a alors un isomorphisme: $$\operatorname{Quot}^G(\lambda, \mu) \simeq \Spec (\Cbb[t]/\langle t^2 \rangle).$$
\end{teo}

\begin{rem}\
\begin{itemize}
 \item Dans le cas où $n=1$, on a $G= \operatorname{SL}(2) \times H $ et $\lambda=(2 \omega_1,0)$.
 \item Dans le cas où $n\geq2$, on a $\lambda=( \omega_1,0)$.
 \end{itemize}
\end{rem}

\subsection{Le schéma $\operatorname{Quot}^G(\lambda, \mu)$ n'a qu'un seul point}

On sait déjà que $\operatorname{Quot}^G(\lambda, \mu)$ admet le point fermé $z$ correspondant à $\Ncal(\lambda, \mu)$. Dans ce paragraphe, on montre qu'il n'y en a pas d'autre.
On commence par quelques rappels sur les algèbres et les modules ``horosphériques''.

\begin{defi}\label{2.2} Soient $R$ une $G$-algèbre, et $V$ un $R$-$G$-module. On dit que $V$ est \textit{horosphérique} si pour tout poids dominants $\lambda_1$, $\lambda_2$, on a $$R_{(\lambda_1)} \cdot V_{(\lambda_2)} \subseteq V_{(\lambda_1+\lambda_2)}.$$
On dit que la $G$-algèbre $R$ est \textit{horosphérique} si elle est horosphérique en tant que $R$-$G$-module. 
\end{defi}

Le théorème suivant découle de \cite{KR} Theorem 3 p 356~:

\begin{teo}\label{2.3} Soit $R$ une $G$-algèbre.
Soit $E \subseteq R$ un sous-$G$-module qui engendre $R$ en tant qu'algèbre.
Alors l'algèbre $R$ est horosphérique si et seulement si pour tout poids dominants $\lambda_1, \lambda_2$ on a  $$E_{(\lambda_1)} \cdot E_{(\lambda_2)} \subseteq R_{(\lambda_1 + \lambda_2)}.$$
\end{teo}

\noindent {\bf Preuve.}
Comme $R$ est la limite inductive de ses sous-algèbres de type fini $G$-stables, il suffit de montrer le théorème dans le cas où l'algèbre $R$ est de type fini.
On peut alors supposer que $E$ est un $G$-module de dimension finie.

Notons $(\Sym E) / I$ le plus grand quotient horosphérique de la $G$-algèbre $\Sym E$.
L'idéal $I$ est homogène, engendré par les composantes isotypiques $$[(\operatorname{S}^mE)_{(\lambda)} \cdot (\operatorname{S}^nE)_{(\mu)}]_{(\nu)},$$ où $m$, $n$ sont des entiers, et $\lambda$, $\mu$, $\nu$ des poids dominants tels que $\lambda + \mu \not = \nu$.

Selon \cite{KR} Theorem 3 p 356, l'idéal $I$ est en fait engendré par sa composante homogène de degré $2$, notée $I_2$.

Ainsi, si on note $J$ le noyau de la surjection canonique $\Sym E \longrightarrow R$, l'algèbre $R$ est horosphérique si et seulement si $I_2 \subseteq J$, c'est-à-dire si et seulement si pour tout poids dominants $\lambda_1, \lambda_2$ on a $E_{(\lambda_1)} \cdot E_{(\lambda_2)} \subseteq R_{(\lambda_1 + \lambda_2)}.$
\hfill $\Box$

\begin{cor}\label{2.4} Soit $R$ une $G$-algèbre horosphérique engendrée par un sous-$G$-module $E \subseteq R$. Soit $V$ un $R$-$G$-module engendré par un sous-$G$-module $W \subseteq V$.
Alors $V$ est un $R$-$G$-module horosphérique si et seulement si pour tout poids dominants $\lambda_1, \lambda_2$ on a  $$E_{(\lambda_1)} \cdot W_{(\lambda_2)} \subseteq W_{(\lambda_1 + \lambda_2)}.$$
\end{cor}

\noindent {\bf Preuve.}
Remarquons que le $R$-$G$-module $V$ est horosphérique si et seulement si la $G$-algèbre $ R \oplus \epsilon V$ (où on pose $\epsilon ^2=0$) est horosphérique. En appliquant le théorème précédent à cette algèbre (engendrée par $E \oplus \epsilon W$), on obtient le corollaire. 
\hfill $\Box$

On établit le lemme suivant à l'aide du corollaire précédent~: 

\begin{lem}\label{2.5} Le $A(\lambda)$-module gradué $N(\lambda, \mu)$ est engendré par sa composante homogène de degré $1$.
\end{lem}

\noindent {\bf Preuve.}
Notons $\langle N(\lambda, \mu)_1 \rangle$ le sous-$A(\lambda)$-module de $N(\lambda, \mu)$ engendré par la composante homogène de degré $1$. Il s'agit de montrer que $\langle N(\lambda, \mu)_1 \rangle = N(\lambda, \mu)$.

On remarque que $A(\lambda)$ est une algèbre horosphérique engendrée par la composante homogène $A(\lambda)_1=V(\lambda)^*$, et que $$\overline{M(\lambda, \mu)} :=M(\lambda, \mu)/\langle N(\lambda, \mu)_1 \rangle$$ est un $A(\lambda)$-module gradué  engendré par sa composante homogène de degré $0$~: $$\overline{M(\lambda, \mu)}_0 \cong V(\mu)^* .$$
Enfin, on a un isomorphisme $$A(\lambda)_1 \cdot  \overline{M(\lambda, \mu)}_0 \cong V(\lambda + \mu)^* .$$
Le module $\overline{M(\lambda, \mu)}$ est donc horosphérique, selon le corollaire précédent.

Si $m$ est un entier, les composantes isotypiques de $A(\lambda)_m \cdot  \overline{M(\lambda, \mu)}_0$ de type différent de $m\lambda^* + \mu^*$ sont donc nulles.
Autrement dit, on a pour tout $m$~: $$N(\lambda, \mu)_m \subseteq \langle N(\lambda, \mu)_1 \rangle,$$ ce qui montre le lemme.
\hfill $\Box$

\begin{prop}\label{2.6}
Le schéma $\operatorname{Quot}^G(\lambda, \mu)$ a un unique point fermé $z$.
\end{prop}
 
\noindent {\bf Preuve.} Soit $P\subseteq M(\lambda, \mu)$ un sous-$A(\lambda)$-module $G$-stable tel qu'on ait un isomorphisme de $G$-modules $$M(\lambda, \mu)/P \simeq \bigoplus_{m=0}^\infty V(m \lambda + \mu)^*.$$
Il s'agit de montrer que $P=N(\lambda, \mu)$.

Si $p$ est un entier, on note $M(\lambda, \mu)_{\geq p}$ le sous-$A(\lambda)$-module gradué $G$-stable $$M(\lambda, \mu)_{\geq p} :=\bigoplus_{m=p}^\infty M(\lambda, \mu)_m \subseteq M(\lambda, \mu).$$ 
Montrons d'abord que $P \subseteq M(\lambda, \mu)_{\geq 1}$.

Par l'absurde, supposons le contraire~: on a alors $$ M(\lambda, \mu)=P + M(\lambda, \mu)_{\geq 1}.$$

Montrons par une récurrence descendante que pour tout $m \in \Nbb$, $$[M(\lambda, \mu)_{\geq m}]_{(\mu^*)} \subseteq P,$$ ce qui donnera une contradiction. 

Si l'entier $m$ est suffisament grand, on a $$[M(\lambda, \mu)_{\geq m}]_{(\mu^*)} = 0,$$ car le $G$-module $M$ est à multiplicités finies.

Soit $m$ un entier tel que $[M(\lambda, \mu)_{\geq m+1}]_{(\mu^*)} \subseteq P$.

Comme $$M(\lambda, \mu)_0 \subseteq P + M(\lambda, \mu)_{\geq 1},$$ on a, en appliquant $A(\lambda)_m$, $$ M(\lambda, \mu)_m =A(\lambda)_m \cdot M(\lambda, \mu)_0 \subseteq A(\lambda)_m \cdot P + M(\lambda, \mu)_{\geq m+1}.$$

Puis, en prenant la composante isotypique de type $\mu^*$~: $$ [M(\lambda, \mu)_m]_{(\mu^*)}  \subseteq A(\lambda)_m \cdot P + [M(\lambda, \mu)_{\geq m+1}]_{(\mu^*)}.$$

D'où le résultat, par récurrence. Ainsi on a $$P \subseteq M(\lambda, \mu)_{\geq 1}.$$

On veut maintenant montrer que $P=N(\lambda, \mu)$~; il suffit pour cela de montrer que $N(\lambda, \mu)$ est inclus dans  $P$.
Selon le lemme précédent, il suffit de montrer que $N(\lambda, \mu)_1$ est inclus dans $P$.

On a vu que les composantes isotypiques non nulles de $N(\lambda, \mu)_1$ sont toutes de type strictement inférieur à $\lambda^*+\mu^*$. Donc leurs images dans $$M(\lambda, \mu)_{\geq 1}/P \simeq \bigoplus_{m=1}^\infty V(m \lambda + \mu)^*$$ sont toutes nulles, et donc $N(\lambda, \mu)_1$ est inclus dans $P$, ce qui montre le lemme.
\hfill $\Box$

\subsection{L'espace tangent au schéma $\operatorname{Quot}^G(\lambda, \mu)$ en $z$}

Le but de ce paragraphe est de démontrer la proposition suivante~:

\begin{prop} \label{et} L'espace tangent en $z$ au schéma Quot
  invariant $\operatorname{Quot}^G(\lambda, \mu)$ est nul, sauf si on
  a (à revêtement fini de $G$ près) 
 $G= \operatorname{Spin}(2n+1) \times H $ pour un groupe réductif
 connexe $H$ et $V(\lambda)=\Cbb^{2n+1}$ et $\mu=(\mu_1, \mu_2)$ avec
 $\langle \mu_1 , \alpha_n^\vee \rangle \geq 1$. 
L'espace tangent est alors de dimension $1$.
\end{prop}

\subsubsection{Une condition nécessaire pour que l'espace tangent soit non nul}

Selon la proposition \ref{1.11}, on a un isomorphisme
$$\operatorname{T}_z \operatorname{Quot}^G(\lambda, \mu) \cong
\operatorname{Hom}_{A(\lambda)}^G(N(\lambda, \mu),Q(\lambda, \mu)).$$

On sait qu'on a une équivalence de catégories abéliennes entre les faisceaux $G$-linéarisés sur l'espace homogène $G \cdot v_\lambda$ et les modules rationnels sur le groupe d'isotropie $G_{v_\lambda}$. Elle est donnée par le foncteur qui à un faisceau $G$-linéarisé $\Fcal$ associe sa fibre $\Fcal_{v_\lambda}$ en $v_\lambda$.
Ainsi, la suite exacte $$0 \longrightarrow \Ncal(\lambda, \mu)|_{G \cdot v_\lambda} \longrightarrow \Mcal(\lambda, \mu)|_{G \cdot v_\lambda} \longrightarrow \Qcal(\lambda, \mu)|_{G \cdot v_\lambda} \longrightarrow 0$$ donne une suite exacte \begin{equation} \label{secfibres} 0 \longrightarrow \Ncal(\lambda, \mu)_{ v_\lambda} \longrightarrow \Mcal(\lambda, \mu)_{ v_\lambda} \longrightarrow \Qcal(\lambda, \mu)_{ v_\lambda} \longrightarrow 0. \end{equation}
De plus, on a un isomorphisme $$ \Hom^G_{\Ocal_{C_\lambda}}(\Ncal(\lambda, \mu)|_{G \cdot v_\lambda}, \Qcal(\lambda, \mu)|_{G \cdot v_\lambda}) \cong \Hom^{G_{ v_\lambda}}(\Ncal(\lambda, \mu)_{ v_\lambda}, \Qcal(\lambda, \mu)_{ v_\lambda}) .$$

\begin{lem}
Le morphisme de restriction\\ $ \Hom^G(\Ncal(\lambda, \mu), \Qcal(\lambda, \mu)) \longrightarrow \Hom^G(\Ncal(\lambda, \mu)|_{G \cdot v_\lambda}, \Qcal(\lambda, \mu)|_{G \cdot v_\lambda})$ est injectif.
\end{lem}

\noindent {\bf Preuve.}
Soit $\phi$ un morphisme non nul de $\Ncal(\lambda, \mu)$ vers  $\Qcal(\lambda, \mu))$ au dessus du cône $C_\lambda$. Comme le cône est affine, il existe une section globale $s$ de $\Ncal(\lambda, \mu)$ telle que $\phi(s) \not = 0$. Puis, comme $Q(\lambda, \mu)$ est un $A(\lambda)$-module sans torsion, la restriction de $\phi(s)$ à tout ouvert non vide de $C_\lambda$ est non nulle. En particulier, $\phi(s|_{G \cdot v_\lambda})=\phi(s)|_{G \cdot v_\lambda}$ est non nulle. 
\hfill $\Box$

\begin{prop} La suite exacte courte de $G_{v_\lambda}$-modules (\ref{secfibres}) s'identifie à la suivante~:
$$0 \longrightarrow (V(\mu)/(V(\mu)^{U_\lambda}))^* \longrightarrow V(\mu)^* \longrightarrow (V(\mu)^{U_\lambda})^*\longrightarrow 0.$$

\end{prop}

\noindent {\bf Preuve.}
La fibre du faisceau $\Mcal = \Ocal_{C_\lambda} \otimes_\Cbb V(\mu)^*$ est $V(\mu)^*$.
Déterminons la fibre de $\Qcal(\lambda, \mu)$ en $v_\lambda$.

On note $f$ et $\pi$ les projections naturelles~: 
$$\xymatrix{    & G / G_{v_\lambda} =G \cdot {v_\lambda} \ar[d]^\pi \\ G/B  \ar[r]^{f~~~~~~~} & G/P_\lambda = G \cdot [v_\lambda] } $$

On remarque que l'on a un isomorphisme d'algèbres $$\Cbb[G//U] \cong \bigoplus_{\nu \in \Lambda^+}\operatorname{H}^0(G/B, \Lcal_\nu),$$ où la multiplication de l'algèbre de droite est celle induite par les multiplications\\ $\Lcal_{\nu_1} \otimes \Lcal_{\nu_2} \longrightarrow \Lcal_{\nu_1+\nu_2}$.

On a donc un isomorphisme de modules 
$$Q(\lambda,\mu) \cong \bigoplus_{m \in \Nbb}
\operatorname{H}^0(G/B, \Lcal_{m\lambda+\mu}).$$ 
Donc 
$$Q(\lambda,\mu) \cong \bigoplus_{m \in \Nbb}
\operatorname{H}^0(G/P_\lambda, \Lcal_{m\lambda} \otimes
f_*\Lcal_\mu),$$ selon la formule de projection. 

La restriction $\Qcal(\lambda, \mu)|_{G \cdot {v_\lambda}}$ est donc
l'image réciproque du faisceau $f_*\Lcal_\mu$ sur $G/P_\lambda$~:
$$\Qcal(\lambda, \mu)|_{G \cdot {v_\lambda}} \cong  \pi^*(f_*\Lcal_\mu).$$ On a donc un isomorphisme sur les fibres~: $$\Qcal(\lambda, \mu)_{v_\lambda} \cong (f_*\Lcal_\mu)_{ [v_\lambda]},$$ avec $$(f_*\Lcal_\mu)_{ [v_\lambda]} \cong  \operatorname{H}^0(P_\lambda/P_\lambda, f_*\Lcal_\mu) \cong  \operatorname{H}^0(P_\lambda/B, \Lcal_\mu).$$

La variété de drapeaux $P_\lambda/B$ est canoniquement isomorphe à $L_\lambda/(B \cap L_\lambda)$. La fibre $(f_*\Lcal_\mu)_{ [v_\lambda]}$ est donc isomorphe à l'espace des sections globales du faisceau $\Lcal_\mu$  sur $L_\lambda/(B \cap L_\lambda)$, donc au $L_\lambda$-module simple de plus grand poids $\mu^*$~: $$(f_*\Lcal_\mu)_{ [v_\lambda]}\cong (V(\mu)^{U_\lambda})^*.$$

Enfin, le premier terme de la suite exacte est donc bien le dual de $V(\mu)/(V(\mu)^{U_\lambda})$.
\hfill $\Box$

L'espace tangent au schéma Quot invariant se plonge donc dans
$$\Hom^{G_{v_\lambda}}((  V(\mu)/(V(\mu)^{U_\lambda}))^*,
(V(\mu)^{U_\lambda})^*) \cong \Hom^{G_{v_\lambda}}(V(\mu)^{U_\lambda},
V(\mu)/(V(\mu)^{U_\lambda})).$$ 
\begin{prop} \label{cn} Si l'espace tangent à
  $\operatorname{Quot}^G(\lambda, \mu)$ en $z$ est non nul, alors à
  revêtement fini de $G$ près, on a un isomorphisme $G \simeq
  \operatorname{Spin}(V) \times H$ pour un groupe réductif connexe $H$,
  et on a $V(\lambda)= V$. 
\end{prop}

\noindent {\bf Preuve.}
On vient de voir que l'espace tangent à $\operatorname{Quot}^G(\lambda, \mu)$ se plonge dans $$E :=\Hom^{G_{v_\lambda}}(V(\mu)^{U_\lambda},  V(\mu)/(V(\mu)^{U_\lambda})).$$
On remarque que comme le tore $T$ agit sur les espaces $V(\mu)^{U_\lambda}$ et $V(\mu)/(V(\mu)^{U_\lambda})$, on a une action de $T$ sur $E$~: si $\phi$ est un élément de $E$ et $t$ un élément de $T$, alors $$t \cdot \phi : v \longmapsto t \cdot \phi ( t ^ {-1} \cdot v ).$$

On suppose que l'espace $E$ est non nul.

Soit $\phi$ un vecteur propre (non nul) de $E$ sous l'action de $T$.
La restriction de cette action au stabilisateur $T_{v_\lambda}$ de $v_\lambda$ dans $T$ est triviale~: le poids de $\phi$ est donc de la forme $d \lambda$, où $d$ est un entier.

Le groupe  $[L_\lambda,L_\lambda]$ dérivé de $L_\lambda$ est contenu dans le groupe d'isotropie $G_{v_\lambda}$, donc le morphisme $\phi$ est $[L_\lambda,L_\lambda]$-équivariant. Comme $V(\mu)^{U_\lambda}$ est un $[L_\lambda,L_\lambda]$-module simple et comme $\phi$ est non nul, $\phi$ est injectif. On a donc $\phi ( v_\mu) \not =0$.

Notons $v$ l'unique antécédent de $\phi ( v_\mu)$ par la projection canonique $V(\mu) \rightarrow V(\mu)/(V(\mu)^{U_\lambda})$ qui soit un $T$-vecteur propre. Son poids est $\mu - d \lambda$.

Le sous-groupe unipotent maximal $U$ de $G$ est contenu dans $G_{v_\lambda}$, donc le morphisme $\phi$ est $U$-équivariant.
Le vecteur $\phi ( v_\mu)$ est donc invariant par $U$, et on a $$\frak{u} \cdot v \subseteq V(\mu)^{U_\lambda}.$$

Comme le vecteur $v$ n'appartient pas à $V(\mu)^{U_\lambda}$, il existe une racine simple $\alpha$ telle que $$e_\alpha \cdot v \not = 0.$$
On a alors $e_\alpha \cdot v \in V(\mu)^{U_\lambda}$.


Une telle racine $\alpha$ ne peut pas être une racine de $L_\lambda$, car sinon le poids de $v$ serait la somme de $\mu$ et de racines de $L_\lambda$, et $v$ appartiendrait à $V(\mu)^{U_\lambda}$.

On remarque donc (comme $v$ est de poids $\mu - d \lambda$), que $d \lambda$ est somme de racines simples de $L_\lambda$ et d'une seule racine simple de $G$ qui n'est pas une racine de $L_\lambda$.

Soit $\alpha$ une racine positive telle que $e_\alpha \cdot v \not = 0$ et que l'on suppose maximale possible.
Montrons que $e_\alpha \cdot v$ est proportionnel à $v_\mu$.

Si ce n'est pas le cas, il existe une racine simple $\beta$ de $L_\lambda$ telle que $$e_\beta \cdot (e_\alpha \cdot v) \not = 0$$ (car $e_\alpha \cdot v$ appartient à $V(\mu)^{U_\lambda}$).
Or $$e_\beta \cdot (e_\alpha \cdot v) = [e_\beta,e_\alpha]\cdot v +  e_\alpha \cdot (e_\beta \cdot v),$$
avec $[e_\beta,e_\alpha]\cdot v=0$ (car on a supposé $\alpha$ maximale possible) et $e_\beta \cdot v=0$ comme on l'a vu~: une contradiction.

~

Les vecteurs $e_\alpha \cdot v$ et $v_\mu$ sont donc proportionnels.
En considérant les poids, on obtient~: $$ \alpha =  d \lambda.$$

Ainsi, $d \lambda$ est une racine positive de $G$ dont l'écriture comme somme de racine simples contient une 
seule racine simple (avec coefficient 1) qui n'est pas une racine de $L_\lambda$.

L'action de $G$ sur $V(\lambda)$ se factorise donc par celle d'un groupe simple (car $\lambda$ est proportionnel à une racine de $G$).

On vérifie facilement que les seuls systèmes de racines simples qui admettent une racine dominante $\alpha$ dont l'écriture comme somme de racine simples contient une seule racine simple (avec coefficient 1) qui n'est pas orthogonale à $\alpha$ sont ceux de type $B_n$, $n\geq 1$, avec $\alpha = \alpha _1+ \alpha _2+...+ \alpha _n $ (c'est-à-dire $\omega_1$ si $n \geq 2$ et $2 \omega_1$ sinon).

Comme un multiple de $\lambda$ doit être une telle racine $\alpha$, il ne reste plus qu'à éliminer le cas où le système de racines est de type $B_1=A_1$ et où $\lambda= \omega_1=1$.

Supposons donc que l'on a $G=\operatorname{SL}(2) \times H$ et  $\lambda=(1,0)$.
La composante homogène de degré $1$ de $N(\lambda, \mu)$ est le noyau du produit de Cartan $V(1,0)\otimes_\Cbb V(\mu_1,\mu_2 ) \longrightarrow V(\mu_1+1,\mu_2 )$. Elle est donc isomorphe à $V(\mu_1-1,\mu_2 )$ si $\mu_1 \not = 0$, nulle sinon. 

Selon le lemme \ref{2.5}, l'espace $N(\lambda, \mu)_1$ engendre le $A(\lambda)$-module $N(\lambda, \mu)$. On a donc une inclusion $$\Hom _{A(1)} ^G(N(\lambda, \mu),Q(\lambda, \mu)) \hookrightarrow \Hom ^G(N(\lambda, \mu)_1 ,Q(\lambda, \mu))=0,$$ donc l'espace tangent est nul dans ce cas.
\hfill $\Box$

\subsubsection{Cas d'un groupe $G$ simple de type $B_n$, $n\geq 1$} \label{2223}

Dans ce paragraphe, on montre la proposition \ref{et} dans le cas où $G=\operatorname{Spin}(2n+1)$ et  $\lambda = \alpha_1 + ... + \alpha_n$, c'est-à-dire $ \omega_1$ si $n\geq 2$ et  $2 \omega_1$ sinon.

\begin{lem} \label{211} La multiplicité de $V(\mu)$ dans la décomposition de $V(\lambda) \otimes_\Cbb V(\mu)$ en somme directe de modules simples est $1$ si $\langle \mu , \alpha_n^\vee \rangle \not= 0$, et $0$ sinon.
\end{lem}

\noindent {\bf Preuve.}
On va utiliser la formule de Weyl (\cite{Bou2} Ch.VIII, \S9), qui donne les poids et leurs multiplicités d'un module $V(\lambda)$ en fonction de $\lambda$.

On note $(e^\nu)_{\nu \in \Lambda}$ la base canonique de l'algèbre $\Zbb[\Lambda]$ du groupe $\Lambda$ à coefficients dans $\Zbb$.

Si $V$ est un $G$-module de dimension finie, on note $m_\nu$ la multiplicité du poids $\nu$ dans $V$ pour tout $\nu \in \Lambda$. Le \textit{caractère} de $V$ est alors $$ \operatorname{ch}(V)  := \sum_{\nu \in \Lambda}m_\nu e^\nu.$$ 

On note $\rho$ la demi-somme des racines positives de $G$. On note enfin $W$ le groupe de Weyl de $G$ relativement au tore maximal $T$, et pour tout élément $w$ de $W$, on note $\epsilon (w)$ la signature de $w$.

Soit $a_\nu$ la multiplicité du module $V(\nu)$ dans la décomposition en somme directe de modules simples de $V(\lambda) \otimes _\Cbb V(\mu)$~: on a $$V(\lambda) \otimes _\Cbb V(\mu) \simeq \bigoplus _ {\nu \in \Lambda^+} a_\nu V(\nu).$$
On a donc, en prenant les caractères~: $$\operatorname{ch}(V(\lambda)) \operatorname{ch}(V(\mu)) = \sum _ {\nu \in \Lambda^+} a_\nu \operatorname{ch}(V(\nu)).$$
Donc, selon la formule de Weyl~: $$\operatorname{ch}(V(\lambda)) \sum _ {w \in W} \epsilon (w)  e^{w(\rho + \mu)}= \sum _ {\nu \in \Lambda^+} \sum _ {w \in W} a_\nu \epsilon (w)e^{w(\rho + \nu)}.$$

On remarque que $a_\mu$ est le coefficient de $e^{\rho + \mu}$ dans chacun des deux membres.

Les poids du module $V(\lambda)$ sont les $\pm \epsilon_i$ où $i=1,...,n$ et le poids nul, chacun avec multiplicité $1$. Le caractère de $V(\lambda)$ est donc $$\operatorname{ch}(V(\lambda)) = e^{\epsilon_1} +...+ e^{\epsilon_n} +1+   e^{-\epsilon_n} +...+  e^{-\epsilon_1}.$$

L'entier $a_\mu$ est donc le cefficient de $e^{\rho + \mu}$ dans l'expression  $$ (e^{\epsilon_1} +...+ e^{\epsilon_n} +1+   e^{-\epsilon_n} + ...+ e^{-\epsilon_1})\sum _ {w \in W} \epsilon (w)  e^{w(\rho + \mu)}.$$

On remarque que $$\pm \epsilon_i + w(\rho + \mu)$$ n'est jamais égal à $\rho + \mu$, sauf si $w$ est la réflexion simple associée à la racine $\alpha_n$ et $\pm \epsilon_i = \epsilon_n$, avec $\langle \mu , \alpha_n^\vee \rangle = 0$.

Le coefficient de $e^{\rho + \mu}$ est donc $1$ si $\langle \mu , \alpha_n^\vee \rangle \not= 0$, et $0$ sinon, d'où le lemme.
\hfill $\Box$

\begin{lem} \label{212} Soit $m \geq 1$ un entier.
Les poids dominants $\nu$ tels que $V(\nu)$ s'injecte dans $V(m \lambda) \otimes_\Cbb V(\mu)$ et qui sont supérieurs ou égaux à $(m-1) \lambda + \mu$ sont de la forme $$ m \lambda + \mu - \alpha_1- \alpha_2-...- \alpha_{i},$$ avec $i=0,...,n$.
\end{lem}

\noindent {\bf Preuve.}
Soit $\nu$ un poids tel que $$(m-1) \lambda + \mu \leq \nu \leq  m \lambda + \mu.$$
Le poids s'écrit donc $$\nu = m \lambda + \mu - \sum _ {j \in J} \alpha_j,$$ où $J$ est une partie de l'intervalle d'entiers $[1;n]$. On note $r$ son cardinal.

Supposons que $V(m \lambda) \otimes_\Cbb V(\mu)$ contient un vecteur primitif de poids $\nu$. On veut montrer que  $J=[1;r]$.

Notons $r_0$ le plus grand entier tel que $J$ contienne $[1;r_0]$. On veut montrer que $r_0=r$. On va raisonner par l'absurde et supposer $r_0<r$.

On rappelle (\cite{Bou2} VIII.7 Exercice 18) que comme $\lambda$ est colinéaire à $\omega_1$, si $j_1, ... , j_s$ sont des entiers de $[1;n]$ distincts deux à deux tels que $$e_{-\alpha_{j_s}} \cdot ...\cdot e_{-\alpha_{j_1}} \cdot v_{m \lambda} $$ est non nul, alors on a $j_i=i$ pour tout $i=1,...,s$.

On pose $$v_s  := e_{-\alpha_{s}} \cdot ...\cdot e_{-\alpha_{1}} \cdot v_{m \lambda}.$$

Soit $v$ un vecteur primitif de $V(m \lambda) \otimes_\Cbb V(\mu)$ de poids $\nu$.
Le vecteur $v$ s'écrit $$v=\sum_{s=1}^{r_0} v_s \otimes w_s,$$ où chaque $w_s$ est un vecteur de $V(\mu)$ de poids $$ \mu - \sum _ {j \in J \setminus [1;s]} \alpha_j.$$

Notons $I$ l'intervalle d'entiers $[1;r_0]$, et $\frak{u}_I$ l'algèbre de Lie du groupe unipotent $U_I$.
En raison de leurs poids respectifs, les vecteurs $v_s$ sont invariants par $U_I$, mais pas les vecteurs $w_s$ (sauf ceux qui sont nuls). Soit un élément $s_0$ de $[1;r_0]$ tel que $w_{s_0}$ est non nul. Soit un élément $x$ de $\frak{u}_I$ tel que $x \cdot w_{s_0}$ est non nul.
On a $$x\cdot v = \sum _ {s=0}^{r_0} v_s \otimes (x \cdot w_{s}) \not = 0,$$ car $(v_0,...,v_{r_0})$ est une famille libre de $V(m   \lambda)$ et $x \cdot w_{s_0}\not = 0$.
Le vecteur $v$ n'est donc pas primitif~: une contradiction.
\hfill $\Box$

\begin{lem} \label{213}
Pour tout entier $m\geq 0$, on note $N^\prime(\lambda, \mu)_m$ l'unique supplémentaire $G$-stable dans $N(\lambda, \mu)_m$ de la composante isotypique de type $(m-1)\lambda+\mu$. 
Alors $$N^\prime(\lambda, \mu) :=\bigoplus_{m \in \Nbb} N^\prime(\lambda, \mu)_m$$ est un sous-$A(\lambda)$-module homogène de $N(\lambda, \mu)$.
\end{lem}

\noindent {\bf Preuve.}
Il suffit de montrer que pour tout entier $m$, on a $$A(\lambda)_1 \cdot N^\prime(\lambda, \mu)_m \subseteq N^\prime(\lambda, \mu)_{m+1}.$$

Pour cela, on va montrer que la composante isotypique de  $$A(\lambda)_1 \otimes _\Cbb  N^\prime(\lambda, \mu)_m$$ de type $m \lambda+\mu$ est nulle.

Soit $\nu$ un poids dominant tel que la composante isotypique de $N^\prime(\lambda, \mu)_m$ de type $\nu$  soit non nulle. Montrons que la composante isotypique de  $$A(\lambda)_1 \otimes _\Cbb V( \nu)$$ de type $m \lambda+\mu$ est nulle.

Si $(m-1)\lambda+\mu$ n'est pas inférieur ou égal à $\nu$, alors $m\lambda+\mu$ n'est pas inférieur ou égal à $\lambda+\nu$, donc  la composante isotypique de  $$A(\lambda)_1 \otimes _\Cbb V( \nu)$$ de type $m \lambda+\mu$ est bien nulle.

Sinon, selon le lemme \ref{212}, on a $$\nu =  m \lambda + \mu -
\alpha_1- \alpha_2-...- \alpha_{r}$$ pour un certain $r \in [1;n-1]$.

Selon \cite{Bou2} VIII.7 Exercice 18, comme $\lambda$ est colinéaire à $\omega_1$, les composantes isotypiques non nulles de $$V(\lambda) \otimes _\Cbb V(m \lambda + \mu - \alpha_1- \alpha_2-...- \alpha_{r})$$ sont soit de type $(m+1) \lambda + \mu - \alpha_1- \alpha_2-...- \alpha_{r}$, soit de type inférieur ou égal à $(m+1) \lambda + \mu - \alpha_1- \alpha_2-...- \alpha_{r}- \alpha_1$. Celle de type $m \lambda+\mu$ est donc nulle.

Ainsi, la composante isotypique de  $$A(\lambda)_1 \otimes _\Cbb  N^\prime(\lambda, \mu)_m$$ de type $m \lambda+\mu$ est bien nulle, et lemme est démontré.
\hfill $\Box$

\begin{lem} \label{214} On suppose $\alpha_n^\vee(\mu) \not= 0$. Pour tout $m \geq 1$, on a un isomorphisme $$N(\lambda, \mu)_m \simeq V((m-1) \lambda + \mu) \oplus N^\prime(\lambda, \mu)_m.$$

\end{lem}

\noindent {\bf Preuve.}
Notons $W_m$ la composante isotypique de $N(\lambda, \mu)_m$ de type $(m-1)\lambda + \mu$.\\
On a $$ N(\lambda, \mu)_m = W_m \oplus N^\prime(\lambda, \mu)_m.$$
Il s'agit de montrer que pour tout $m$, le $G$-module $W_m$ est simple.

Selon le lemme \ref{2.5}, on a $$N(\lambda, \mu)_m= A(\lambda)_{m-1} \cdot N(\lambda, \mu)_1.$$
Selon le lemme \ref{213}, on a $$A(\lambda)_{m-1} \cdot N^\prime(\lambda, \mu)_1 \subseteq N^\prime(\lambda, \mu)_{m}.$$

On a donc une surjection de $G$-modules $$A(\lambda)_{m-1} \otimes _\Cbb W_1 \twoheadrightarrow W_m.$$

Les $G$-modules $A(\lambda)_{m-1}$ et $W_1$ sont isomorphes respectivement à $V((m-1)\lambda)$ et $V(\mu)$ (selon le lemme \ref{211}), donc $W_m$ est soit nul soit isomorphe à $V((m-1)\lambda + \mu)$.

Le $A(\lambda)$-module $M(\lambda, \mu)$ est sans torsion (car libre).
Dons si $a$ est un vecteur primitif de $A(\lambda)_1$ et $v$ un vecteur primitif de $W_1$, le vecteur $av$ est un vecteur primitif (non nul) de $N(\lambda, \mu)_m$ de poids $(m-1)\lambda + \mu$, donc $W_m$ est en fait isomorphe à $V((m-1)\lambda + \mu)$.
\hfill $\Box$

~

Vérifions maintenant que la proposition \ref{et} est vraie dans notre situation.\\

Le $A(\lambda)$-module  $N(\lambda, \mu)$ est engendré par sa composante homogène de degré $1$.
On a donc un plongement $$\Hom _{A(\lambda)} ^G(N(\lambda, \mu),Q(\lambda, \mu)) \hookrightarrow \Hom ^G(N(\lambda, \mu)_1 ,Q(\lambda, \mu)).$$

Selon le lemme \ref{211}, l'espace $\Hom ^G(N(\lambda, \mu)_1 ,Q(\lambda, \mu))$ est de dimension $1$ si 
$\alpha_n^\vee(\mu) \not= 0$, et nul sinon.

Le théorème est donc vérifié si $\alpha_n^\vee(\mu)=0 $.
Il reste à vérifier que l'espace tangent au schéma Quot invariant est non nul si $\alpha_n^\vee(\mu) \not= 0$.

Selon le lemme \ref{214}, le $A(\lambda)$-module $N(\lambda, \mu)/N^\prime(\lambda, \mu)$ admet alors $h_{\lambda, \mu}$ comme fonction de Hilbert.
De plus il engendré par sa composante homogène de degré $1$~: c'est donc un quotient de $A(\lambda)\otimes V(\mu)=M(\lambda, \mu)$. Il est donc isomorphe à $Q(\lambda, \mu)$, selon la proposition \ref{2.6}. 

Donc si $\alpha_n^\vee(\mu)$ est non nul, il existe un morphisme non nul de $N(\lambda, \mu)$ vers $Q(\lambda, \mu)$, et l'espace tangent est donc non nul~: la proposition \ref{et} est vraie dans notre situation.

\subsubsection{Cas général}

 \label{casg}

On va conclure à l'aide du fait général suivant:

\begin{lem} \label{produit}
Soient $G_1$ et $G_2$ deux groupes réductifs connexes, chacun muni d'un tore maximal et d'un sous-groupe de Borel le contenant. On note $\Lambda _i ^+$ l'ensemble des poids dominants de $G_i$.

Soit $X$ un $G_1$-schéma affine,  $\Mcal_1$ un $\Ocal_X$-module cohérent et $G_1$-linéarisé, et $h_1$ une fonction sur  $\Lambda _1 ^+$ à valeurs entières.

Soit $\mu_2$ un poids dominant de $G_2$. On note $V_{G_2}(\mu_2)$ le $G_2$-module simple associé, et $ \Lambda _1 ^+ \times \Lambda _2 ^+ \longrightarrow \Nbb$ la fonction donnée par $h(\lambda_1,\lambda_2):=  h_1(\lambda_1)$ si $\lambda_2 = \mu_2$ et $0$ sinon.

On a un isomorphisme canonique: $$ \operatorname{Quot}_h^{G_1 \times G_2} ( X, \Mcal_1 \otimes_\Cbb V_{G_2}(\mu_2)) \cong \operatorname{Quot}_{h_1}^{G_1} ( X, \Mcal_1).$$
\end{lem}

\noindent {\bf Preuve.} 
On a un morphisme fonctoriel $$  \Quot_h^{G_1 \times G_2} ( X, \Mcal_1 \otimes_\Cbb V_{G_2}(\mu_2)) \longrightarrow  \Quot_{h_1}^{G_1} ( X, \Mcal_1)$$ donné pour tout schéma $S$ par l'application

 $$\begin{array}{rlc}  \Quot_h^{G_1 \times G_2} ( X, \Mcal_1 \otimes_\Cbb V_{G_2}(\mu_2))(S)& \longrightarrow& \Quot_{h_1}^{G_1} ( X, \Mcal_1)  (S)     \\ \Ncal_1  &\longmapsto & \Ncal_1 \otimes _\Cbb V_{G_2}(\mu_2)
\end{array}$$

Cette application est bijective, car tout sous-$G_1 \times G_2$-module de $\Mcal_1 \otimes _\Cbb V_{G_2}(\mu_2)$ est de la forme $\Ncal_1 \otimes _\Cbb V_{G_2}(\mu_2)$, où $\Ncal_1$ est un sous-$G_1$-module de $\Mcal_1$.
\hfill $\Box$\\

Montrons maintenant la proposition \ref{et}.

Selon la proposition \ref{cn}, on peut supposer que l'on a $G = \operatorname{Spin}(2n+1) \times H$, et $\lambda$ est $ ( \omega_1,0)$ si $n \geq 2$ et $(2\omega_1,0)$ si $n=1$. 

On remarque que l'on a, en notant $V_H(\mu_2)$ le $H$-module simple de plus grand poids $\mu_2$~:
$$M(\lambda, \mu) = M(\lambda, (\mu_1, 0)) \otimes_\Cbb V_H(\mu_2).$$
L'action de $H$ sur $X$ est triviale.
On conclut en appliquant le lemme précédent.

\subsection{Détermination de $\operatorname{Quot}^G(\lambda,\mu)$}

On va maintenant établir le théorème \ref{quot}.
Selon les propositions \ref{2.6} et \ref{et}, le schéma $\operatorname{Quot}^G(\lambda,\mu)$ est soit un point réduit, soit isomorphe à $\Spec \Cbb[t]/\langle t ^ n \rangle$, pour un entier $n$ supérieur ou égal à $2$.
Il ne reste plus qu'à montrer que cet entier $n$ est toujours égal à $2$.\\

Notons $t$ une indéterminée, et $\epsilon$ (resp. $\delta$) sa classe dans l'algèbre  $\Cbb[t]/\langle t ^ 2 \rangle$ (resp. $\Cbb[t]/\langle t ^ 3 \rangle$). On identifie $\Spec \Cbb[\epsilon]$ à un sous-schéma fermé de $\Spec \Cbb[\delta]$ à l'aide de la surjection canonique $\Cbb[\delta] \longrightarrow \Cbb[\epsilon]$.

Il suffit de montrer que la différentielle (entre les espaces tangents de Zariski) de tout morphisme $$\Psi :\Spec \Cbb[\delta] \longrightarrow \operatorname{Quot}^G(\lambda,\mu)$$ est nulle, c'est-à-dire que la restriction de $\Psi$ à $\Spec \Cbb[\epsilon]$ correspond au vecteur tangent nul en $z$.

Soit $$\Phi : \Spec \Cbb[\epsilon] \longrightarrow \operatorname{Quot}^G(\lambda,\mu)$$ un morphisme.
La proposition suivante décrit l'ensemble des morphismes qui prolongent $\Phi$ à $\Spec \Cbb[\delta]$.
On reprend les notations de la proposition \ref{1.11} et sa démonstration: le morphisme $\Phi$ correspond à un sous-module $L$ de $\Cbb[\epsilon]\otimes_\Cbb M(\lambda,\mu)$, qui est donné par un morphisme noté $\phi: N(\lambda,\mu) \longrightarrow Q(\lambda,\mu)$ à l'aide de l'expression (\ref{formule}).

\begin{prop}\label{delta}
L'ensemble des morphismes $\Spec \Cbb[\delta] \xrightarrow{~~\Psi~~} \operatorname{Quot}^G(\lambda,\mu)$ dont la restriction à  $\Spec \Cbb[\epsilon]$ est $\Phi$ est en bijection avec l'ensemble des morphismes de $A(\lambda)$-$G$-modules $$\psi : L \longrightarrow  Q(\lambda,\mu)$$ dont la restriction à ${\epsilon L \cong N(\lambda, \mu)}$ coïncide avec $\phi$.
\end{prop}

\noindent {\bf Preuve.}
Le raisonnement est analogue à celui de la preuve de la démonstration \ref{1.11}.

Un morphisme $\Phi$ comme dans l'énoncé correspond à un sous-$A(\lambda)[\delta]$-$G$-module $$K \subseteq \Cbb[\delta] \otimes_\Cbb M(\lambda,\mu) = M(\lambda,\mu) \oplus \delta M(\lambda,\mu) \oplus \delta^2 M(\lambda,\mu)$$ tel qu'on ait l'identification $$\Cbb[\epsilon] \otimes_{\Cbb[\delta]} K \cong L$$ et que le quotient $$ (\Cbb[\delta] \otimes_\Cbb M(\lambda,\mu)) / K $$ soit un $\Cbb[\delta]$-module plat.

Soit $K$ un tel sous-module.

On note $L^\prime$ le sous-$A(\lambda)$-$G$-module (isomorphe à $L$) : $$L^\prime:= \bigcup_{n \in N(\lambda,\mu)} (n + \delta \phi(n)) \subseteq M(\lambda,\mu) \oplus \delta M(\lambda,\mu).$$

Comme précédemment, la première condition donne
\begin{equation}\label{pullback2}
(K+\delta^2 M(\lambda,\mu)) \cap (M(\lambda,\mu) \oplus \delta M(\lambda,\mu)) = L^\prime.      
\end{equation}
Notons que cette condition implique que $\delta^2 K = \delta^2 N(\lambda,\mu)$.

La seconde condition donne
\begin{equation}\label{platitude2}
 K \cap \delta^2 M(\lambda,\mu) = \delta^2 N(\lambda,\mu)
\end{equation}
(en utilisant cette fois le fait que dans un $\Cbb[\delta]$-module plat, les éléments annulés par $\delta^2$ sont ceux ``multiples'' de $\delta$).

Donc pour tout élément $l$ de $L^\prime$, il existe un unique élément $\psi(l)$ de $M(\lambda,\mu)/N(\lambda,\mu)$ (on voit cet élément comme une partie de $M(\lambda,\mu)$) tel que $$l+\delta^2 \psi(l) \subseteq K.$$
On a alors
 \begin{equation}\label{formule2}K= \bigcup_{l \in L^\prime} (l+\delta^2 \psi(l)).\end{equation}
Comme $K$ est un $A(\lambda)$-$G$-module, l'application $\psi$ est un morphisme de $A(\lambda)$-$G$-modules.

Enfin, comme $K$ est stable par multiplication par $\delta$, il contient $\delta L^\prime$. Cela signifie que l'application $\psi|_{\delta N}$ coïncide avec $\phi$ (c'est-à-dire que pour tout $n$, $\psi (\delta n)=\phi (n)$).

Réciproquement, un tel morphisme $\psi$  définit bien via l'expression (\ref{formule2}) un morphisme de $\Spec \Cbb[\delta]$ dans $\operatorname{Quot}  ^G (\lambda,\mu)$. 
\hfill $\Box$\\

Dans la suite, on montre que lorsque le morphisme $\phi$ est non nul, il n'existe pas de tel $\psi$.

\subsubsection{Cas d'un groupe $G$ simple de type $B_n$, $n\geq 1$} 

Dans ce paragraphe, on montre le théorème \ref{quot} dans le cas où $G=\operatorname{Spin}(2n+1)$ et  $\lambda = \alpha_1 + ... + \alpha_n$, c'est-à-dire $ \omega_1$ si $n\geq 2$ et  $2 \omega_1$ sinon.

On pose, pour tout $1 \leq i \leq n$: $$\mu_i := \langle \mu , \alpha_i^\vee \rangle,$$ de sorte qu'on a $\mu = \mu_1 \omega_1 + ... + \mu_n \omega_n$.

On se place dans la situation où le schéma $\operatorname{Quot}  ^G (\lambda,\mu)$ n'est pas un point réduit. On a donc $\mu_1 \geq 0$, ... , $\mu_{n-1} \geq 0$ et  $\mu_n \geq 1$. 

On note $r(\mu)$ le nombre de ces inégalités qui sont en fait des égalités.\\

Les deux lemmes suivants se démontrent à l'aide de la formule de Weyl exactement comme le lemme \ref{211}.
On rappelle qu'on a la décomposition : $\operatorname{S}^2 V(\lambda) = V( 2 \lambda) \oplus V(0)$, de sorte que le caractère de $V( 2 \lambda)$ est $$\operatorname{ch}(V( 2 \lambda))=n + \sum _ i (e^{\epsilon_i } +e^{-\epsilon_i }) +  \sum _{i \not = j}  e^{\epsilon_i- \epsilon_j }    +              \sum _ {i \leq j} (e^{\epsilon_i+ \epsilon_j } +  e^{-\epsilon_i- \epsilon_j }).$$

\begin{lem} \label{q1} Soit $\nu$ un poids dominant, et $i$ un entier tel que $1 \leq i \leq n$.

On suppose que $\nu + \epsilon_i$ (resp. $\nu - \epsilon_i$) est un poids dominant.

Alors la multiplicité de $V(\nu + \epsilon_i)$ (resp. $V(\nu - \epsilon_i)$) dans la décomposition de $V(\lambda)  \otimes  V(\nu)$ en somme directe de modules simples est $1$.
\end{lem}

\begin{lem} \label{q2} La multiplicité de $V(\mu)$ dans la décomposition de $V(2 \lambda) \otimes V(\mu)$ en somme directe de modules simples est $n-r(\mu)$.
\end{lem}

On aura aussi besoin du lemme suivant:

\begin{lem} \label{q3}
Soit $\nu$ un poids dominant.
Soit $\sigma$ une somme de racines positives telle que le $G$-module $V(\lambda)  \otimes  V(\nu)$ contienne un $B$-vecteur propre $v$ de poids $\lambda + \nu - \sigma$.

Soit une écriture du vecteur $v$ de la forme $$ v = \sum _ {\sigma_1 + \sigma_2 = \sigma} w_{\lambda -\sigma_1 } \otimes w^\prime_{\nu -\sigma_2 }$$ où les $\sigma_1$ et les  $\sigma_2$ sont des sommes de racines positives, et chaque vecteur $w_{\lambda -\sigma_1 }$ (resp. $ w^\prime_{\nu -\sigma_2 }$) est un vecteur de $V(\lambda)$ (resp. $ V(\nu)$) de poids $\lambda -\sigma_1$ (resp. $ \nu -\sigma_2 $).

Alors le terme $ w_{\lambda -\sigma } \otimes w^\prime_{\nu}$ est non nul.
\end{lem}

\noindent {\bf Preuve.} Soit $\overline{\sigma_2}$ une somme de racines positives telle que le terme $w_{\lambda -\overline{\sigma_1} } \otimes w^\prime_{\nu -\overline{\sigma_2} }$ soit non nul  et minimale pour cette propriété (en posant $\overline{\sigma_1}:= \sigma - \overline{\sigma_2}$).

On suppose par l'absurde que $\overline{\sigma_2}$ est non nulle.

Il existe alors une racine simple $\alpha$ telle $e_\alpha \cdot w^\prime_{\nu -\overline{\sigma_2 }}$ soit non nul.

Or $e_\alpha \cdot v $ est nul, donc on a $$  \sum _ {\sigma_1 + \sigma_2 = \sigma} [(e_\alpha \cdot w_{\lambda -\sigma_1 }) \otimes w^\prime_{\nu -\sigma_2 }    +  w_{\lambda -\sigma_1 } \otimes (e_\alpha \cdot w^\prime_{\nu -\sigma_2 } )]  =0.$$

Donc le terme $w_{\lambda -\overline{\sigma_1} } \otimes (e_\alpha \cdot w^\prime_{\nu -\overline{\sigma_2 }} )$ est nul (d'après la décomposition en somme directe $V(\lambda)  \otimes  V(\nu) = \bigoplus _{\sigma_1,~\sigma_2}V(\lambda)_{\sigma_1}  \otimes  V(\nu)_{\sigma_2}$) : une contradiction.
\hfill $\Box$

\begin{prop} \label{q4} La composante isotypique de type $\mu$ de $M(\lambda,\mu)_2$ est incluse dans le $A(\lambda)$-module engendré par $N^\prime(\lambda, \mu)_1$ (cette notation a été introduite dans le lemme \ref{214}).
\end{prop}

\noindent {\bf Preuve.} Selon le lemme \ref{q2}, il suffit de montrer que l'espace $A(\lambda)_1 \cdot N^\prime(\lambda, \mu)_1$ contient $n-r(\mu)$ vecteurs propres de poids $\mu$ pour l'action de $B$ qui sont linéairement indépendants.

On identifie $\frak{g}$ à l'algèbre de Lie des matrices de taille $(2n+1) \times (2n+1) $ antisymétriques par rapport à la seconde diagonale.
On identifie alors $V(\lambda) $ à l'espace vectoriel $\Cbb ^ {2n+1}$, dont la base canonique est notée $(e_1, ... , e_n, e_0, e_{-n}, ... , e_{-1})$.

On considère le module simple $V(2 \lambda)$  comme le quotient  de $V(\lambda) \otimes V(\lambda)$ par son sous-espace vectoriel $W$ engendré par $ e_1 \otimes e_{-1}+ ... + e_n  \otimes e_{-n} + \frac{1 }{2} e_0 \otimes e_0 $ et par la famille $(e_i \otimes e_j - e_j \otimes e_i)_{i,j}$: $$V(2 \lambda) = (V(\lambda) \otimes V(\lambda))/W.$$

Notons $I$ l'ensemble des entiers $i$ (avec $ 1 \leq i \leq n$) tels que $\mu - \epsilon_i$ soit un poids dominant.
On remarque que \begin{itemize}
 \item Un entier $i$ différent de $n$ appartient à $I$ si et seulement si $\mu_i \geq 1$.
 \item L'entier $n$ appartient à $I$ si et seulement si $\mu_n \geq 2$.
 \end{itemize}
Le cardinal de $I$ est donc: $$\operatorname{card}(I) = n - r(\mu).$$

On va associer à tout élément de $I$ un $B$-vecteur propre (que l'on notera $\overline{v_i^\prime}$) de poids $\mu$ dans $M(\lambda, \mu)_2$.\\

Soit $i$ un élément de $I$.
Selon le lemme \ref{q1}, il existe un $B$-vecteur  propre $v_i$ de poids $\mu - \epsilon_i$ dans $M(\lambda, \mu)_1 = V(\lambda) \otimes V(\mu)$ (donc en fait dans $N^\prime(\lambda, \mu)_1$). Selon le lemme \ref{q3}, celui-ci s'écrit (à un scalaire non nul près) 
\begin{equation}\label{equa1}
v_i = e_{-i} \otimes v_\mu +         e_{-(i+1)} \otimes w^i_{-(i+1)} +...+        e_{-n} \otimes w^i _{-n}+        e_{0} \otimes w^i _0     +   e_{n} \otimes w^i_{n}   +...+        e_{1} \otimes w^i_1 
\end{equation}
où les $w_j^i$ sont des $T$-vecteurs propres de $V(\mu)$ de poids strictement inférieurs à $\mu$. (On rappelle que $v_\mu$ est un $B$-vecteur propre de $V(\mu)$.)

Posons $\nu := \mu - \epsilon _ i$.
Selon le lemme \ref{q1}, le $G$-module $V(\mu)$ s'injecte dans $V(\lambda) \otimes V(\nu)$.

En appliquant à nouveau le lemme \ref{q3}, on obtient un $B$-vecteur propre de poids $\mu$ dans $V(\lambda) \otimes V(\lambda) \otimes  V(\mu)$ qui s'écrit
\begin{equation}\label{equa2}
v_i^\prime = e_{i} \otimes v_i +         e_{i-1} \otimes (u^i_{-\alpha_{i-1}} v_i) +...+        e_{1} \otimes( u^i_{-\alpha_{1}-...-\alpha_{i-1}} v_i)
\end{equation}
où chaque $u^i_{\nu^\prime}$ est un vecteur propre de poids $\nu^\prime$ pour l'action de $T$ dans l'algèbre enveloppante de $\frak{u}^- :=\bigoplus_{\alpha \in R_+} \frak{g}_{-\alpha}$.\\

On note $\overline{v_i^\prime}$ la classe de $v_i^\prime$ modulo $W \otimes V(\mu)$ : on voit donc $\overline{v_i^\prime}$ comme un élément de $M(\lambda, \mu)_2 = V(2 \lambda) \otimes V(\mu)$.

On remarque que, par construction, le vecteur $\overline{v_i^\prime}$ est un $B$-vecteur propre de poids $\mu$ et appartient à $A(\lambda)_1 \cdot N^\prime(\lambda, \mu)_1$.

Il ne reste plus qu'à montrer que la famille $(\overline{v_i^\prime})_{i \in I}$ ainsi construite est libre.

Pour cela, montrons que les $v_i$ sont en fait linéairement indépendants modulo le sous-espace $$W^\prime :=  W \otimes V(\mu)  + V(\lambda) \otimes V(\lambda) \otimes  (V(\mu)_{<\mu}),$$ où $V(\mu)_{<\mu}$ est le sous-espace vectoriel de $V(\mu)$ engendré par ses $T$-vecteurs propres de poids strictement inférieurs à $\mu$.

En remplaçant $v_i$ dans (\ref{equa2}) par son expression  (\ref{equa1}), on obtient que $v_i^\prime$ est congru modulo $W^\prime$ à $$  e_i \otimes e_{-i} \otimes v_\mu + x_{i-1} e_{i-1} \otimes e_{-(i-1)} \otimes v_\mu + ... + x_1 e_1 \otimes e_{-1} \otimes v_\mu,$$ où les $x_j$ sont des scalaires.

D'où le résultat.
\hfill $\Box$\\

On peut maintenant montrer que le théorème \ref{quot} est vérifié dans notre situation.
Pour cela, avec les notations de la proposition \ref{delta}, on se donne un morphisme $\phi$  non nul, et on montre qu'il n'existe aucun morphisme $\psi$ tel que dans la proposition. 

Soient $a_1, ...~ ,a_s$ des éléments de $A(\lambda)_1=V(\lambda)$ et $v_1, ...~ ,v_s$ des éléments de la composante isotypique de $M(\lambda, \mu)_1$ de type $\mu$ tels que le vecteur $$\sum_j a_j \otimes v_j$$ soit un $B$-vecteur propre de poids $\mu$ (selon le lemme \ref{211}, de tels éléments existent, car $\mu_n$ est non nul).  

Pour tout $j$, on note $v_j ^\prime$ l'unique représentant dans $M(\lambda, \mu)_0$ de $\phi(v_j)$, et on pose $$l_j :=  v_j + \epsilon v_j ^\prime.$$
C'est un élément de $L$. On pose enfin $$l:= \sum_j a_j \cdot l_j.$$

Le vecteur $l$ appartient à $A(\lambda)_1 \cdot L$, et est soit nul, soit un $B$-vecteur propre de poids $\mu$ (par construction).

Son image par $\psi$ appartient donc à $A(\lambda)_1 \cdot Q(\lambda, \mu) $, et est donc nulle (car ce dernier espace ne contient pas de $B$-vecteur propre de poids $\mu$).\\

Le vecteur $l$ se décompose sous la forme $$ l =  l ^{\prime } + l ^{\prime \prime} ,$$
en posant $$l^{\prime }:= \sum_j a_j \cdot v_j$$ et $$l ^{\prime \prime}:= \epsilon \sum_j a_j \cdot v_j ^{\prime }.$$

Le vecteur $l^{\prime }$ appartient à $M(\lambda, \mu)_2$, et est soit nul, soit un $B$-vecteur propre de poids $\mu$.
Selon la proposition \ref{q4}, il appartient donc à $A(\lambda)_1 \cdot L $ (car $N^{\prime }(\lambda, \mu)$ est inclus dans $L$), et son image par $\psi$ est donc nulle (de même que celle de $l$).\\

Le vecteur $l ^{\prime \prime}$ appartient donc lui aussi à $L$.
En fait, c'est un $B$-vecteur propre de poids $\mu$ (non nul, car le morphisme naturel $A(\lambda)_1 \otimes M(\lambda, \mu)_0 \longrightarrow A(\lambda)_1 \cdot M(\lambda, \mu)_0 $ est en fait un isomorphisme)  appartenant à $\epsilon N(\lambda, \mu)_1$.

Son image par $\psi$ est donc non nulle (car la restriction de $\psi$ à $\epsilon N(\lambda, \mu) $ coïncide avec $\phi$).

On obtient ainsi une contradiction avec $$\psi(l ^{\prime \prime})=  \psi(l ) - \psi(l ^\prime)=0.$$
Il n'existe donc pas de tel $\psi$, et le théorème est vérifié.

\subsubsection{Conclusion}

On en déduit le théorème à l'aide du lemme \ref{produit}, comme dans le \S\ref{casg}.

\end{document}